\documentclass[aop,MSNbibl,seceqn,citesort,dvips]{arximspdf}

\doi{10.1214/11-AOP696}
\volume{41}
\issue{3A}
\pubyear{2013}
\firstpage{1315}
\lastpage{1361}

\makeatletter

\newcommand{\eps}{\varepsilon}

\renewcommand{\P}{\mathcal{P}}
\newcommand{\E}{\mathcal{E}}
\newcommand{\M}{\mathcal{M}}
\newcommand{\ch}{\operatorname{ch}}
\renewcommand{\th}{\operatorname{th}}
\newcommand{\Av}{\operatorname{Av}}

\newcommand{\Reals}{\mathbb{R}}
\newcommand{\Natural}{\mathbb{N}}
\newcommand{\vsi}{{\bolds{\sigma}}}
\newcommand{\vrho}{{\bolds{\rho}}}

\newtheorem{lemma}{Lemma}
\newtheorem{theorem}{Theorem}

\newproclaim{Remark}{Remark}
\newproclaim{Example}{Example}

\newcommand{\la}{\langle}
\newcommand{\LL}{\mathcal{L}}

\newcommand{\ra}{\rangle}

\makeatother

\begin{document}
\begin{frontmatter}

\title{Spin glass models from the point of view of spin distributions}
\runtitle{Spin distributions in spin glasses}

\begin{aug}
\author[A]{\fnms{Dmitry} \snm{Panchenko}\corref{}\thanksref{t1}\ead[label=e1]{panchenk@math.tamu.edu}}
\runauthor{D. Panchenko}
\affiliation{Texas A\&M University}
\address[A]{Department of Mathematics\\
Texas A\&M University\\
Mailstop 3386, Room 209\\
College Station, Texas 77843\\
USA\\
\printead{e1}} %adresu isvedimo komanda gale!
\end{aug}

\thankstext{t1}{Supported in part by an NSF grant.}

% HISTORY:
\received{\smonth{11} \syear{2010}}
\revised{\smonth{6} \syear{2011}}

% ABSTRACT
%
\begin{abstract}
In many spin glass models, due to the symmetry among sites, any
limiting joint distribution of spins under the annealed Gibbs measure
admits the Aldous--Hoover representation encoded by a function $\sigma:
[0,1]^4\to\{-1,+1\}$, and one can think of this function as a generic
functional order parameter of the model. In a class of diluted models,
and in the Sherrington--Kirkpatrick model, we introduce novel
perturbations of the Hamiltonian that yield certain invariance and
self-consistency equations for this generic functional order parameter
and we use these invariance properties to obtain representations for
the free energy in terms of $\sigma$. In the setting of the
Sherrington--Kirkpatrick model, the self-consistency equations imply
that the joint \mbox{distribution} of spins is determined by the joint
distributions of the overlaps, and we give an explicit formula for
$\sigma$ under the Parisi ultrametricity hypothesis. In addition, we
discuss some connections with the Ghirlanda--Guerra identities and
stochastic stability and describe the expected Parisi ansatz in the
diluted models in terms of $\sigma$.
\end{abstract}

% KEYWORDS
%

\end{frontmatter}

%s1 ###
\section{Introduction and main results}

In various mean-field spin glass models, such as the
Sherrington--Kirkpatrick model and
diluted $p$-spin and $p$-sat models that we will focus on in this
paper, one considers a
random Hamiltonian $H_N(\vsi)$ indexed by spin configurations $\vsi
\in\Sigma_N = \{-1,+1\}^N$ and defines the corresponding Gibbs
measure $G_N$ as a random probability measure on $\Sigma_N$ given by
%
%e1.1 ###
\begin{equation}
G_N(\vsi) = \frac{1}{Z_N} \exp(- H_N(\vsi)),
\end{equation}
where the normalizing factor $Z_N$ is called the partition function.
Let $(\vsi^l)_{l\geq1}$
be an i.i.d. sequence of replicas from measure $G_N$. Let $\mu_N$
denote the joint distribution
of the array of all spins on all replicas, $(\sigma_i^l)_{1\leq i\leq
N, 1\leq l}$, under the
annealed product Gibbs measure $\mathbb{E}G_N^{\otimes\infty}$ which
means that for any choice
of signs $a_i^l\in\{-1,+1\}$, and for any $n\geq1$,
%
%e1.2 ###
\begin{eqnarray}\label{muN}
&&
\mu_N(\{\sigma_i^l = a_i^l\dvtx 1\leq i\leq N, 1\leq l\leq
n\})\nonumber\\[-8pt]\\[-8pt]
&&\qquad=
\mathbb{E}G_N^{\otimes n}(\{\sigma_i^l = a_i^l\dvtx 1\leq i\leq N, 1\leq
l\leq n\}).\nonumber
\end{eqnarray}
In most mean-field spin glass models this distribution has the
following two symmetries. Clearly, it is always invariant under the
permutation of finitely many replica indices $l\geq1$, but in most
models $\mu_N$ is also invariant under the permutation of coordinates
$i\in\{1,\ldots,N\}$ since the distribution of $H_N(\vsi)$ is symmetric
under the permutation of coordinates of $\vsi$, and this invariance of
$\mu_N$ is called symmetry among sites. Let us think of $\mu_N$ as a
distribution on $(\sigma_i^l)$ for all $i,l\geq1$ simply by setting
$\sigma_i^l=0$ for $i>N$. It is usually not known how to prove that the
sequence $(\mu_N)$ converges (in the sense of convergence of
finite-dimensional distributions) and, in fact, even the answer to a
much less general question whether the distribution of one overlap
$N^{-1}\sum _{i\leq N}\sigma_i^1\sigma_i^2$ under
$\mathbb{E}G_N^{\otimes2}$ converges is known only in the
Sherrington--Kirkpatrick model with all $p$-spin interaction terms
present, the proof of which relies on the Parisi formula for the free
energy; see~\cite{T-P,PM}. As a result, we will consider a family $\M$
of all possible limits over the subsequences of~$(\mu_N)$. Whenever we
have symmetry among sites, any limiting distribution $\mu\in\M$ will be
invariant under the permutations of both row and column coordinates $l$
and~$i$. Such two-dimensional arrays are called exchangeable arrays and
the representation result of Aldous~\cite{Aldous} and Hoover
\cite{Hoover} (see also~\cite{Austin}) states that there exists a
measurable function $\sigma_\mu\dvtx[0,1]^4\to\Reals$ such that the
distribution $\mu$ coincides with the distribution of the array
$(s_i^l)$ given by
%
%e1.3 ###
\begin{equation}\label{sigma}
s_i^l=\sigma_\mu(w,u_l,v_i,x_{i,l}),
\end{equation}
where random variables $w,(u_l), (v_i), (x_{i,l})$ are i.i.d. uniform
on $[0,1]$.
This function $\sigma_\mu$ is defined uniquely up to some
measure-preserving transformations
(Theorem 2.1 in~\cite{Kallenberg}) so we can identify the distribution
$\mu$ of array $(s_i^l)$ with
the function $\sigma_\mu$. Since we only consider the case when spins
and thus $\sigma_\mu$
take values in $\{-1,+1\}$, the distribution $\mu$ is completely
encoded by the function
%
%e1.4 ###
\begin{equation}\label{fop}
\bar{\sigma}_\mu(w,u,v) = \mathbb{E}_x \sigma_\mu(w,u,v,x),
\end{equation}
where $\mathbb{E}_x$ is the expectation in $x$ only and we can think
of this last coordinate as a dummy
variable that generates a Bernoulli r.v. with expectation $\bar{\sigma
}_\mu(w,u,v)$. However,
keeping in mind that a function of three variables $\bar{\sigma}_\mu
$ encodes the distribution
of the array (\ref{sigma}), for convenience of notation we will
sometimes not identify a Bernoulli
distribution with its expectation (especially, in the diluted models)
and work with the function
$\sigma_\mu(w,u,v,x)$.

One can think of a function $\sigma_\mu$ (or $\bar{\sigma}_\mu$)
as what physicists might
call a generic ``functional order parameter'' of the model, and it is
easy to see that information
encoded by $\sigma_\mu$ is equivalent to the limiting joint
distribution of all multi-overlaps
%
%e1.5 ###
\begin{equation}\label{multioverlapN}
R_{l_1,\ldots, l_n}^N = N^{-1} \sum_{1\leq i\leq N} \sigma
_i^{l_1}\cdots\sigma_i^{l_n}
\end{equation}
for all $n\geq1$ and all $l_1,\ldots, l_n\geq1$ under $\mu_N$,
which may be a more
familiar object than the joint distribution of spins. Indeed, by
expanding the powers of
(\ref{multioverlapN}) in terms of products of spins and using symmetry
among sites,
in the limit one can express the joint moments of multi-overlaps in
terms of the joint
moments of spins and vice versa. By comparing these moments, the
asymptotic joint
distribution of (\ref{multioverlapN}) over a subsequence of $\mu_N$
converging to
$\mu$ coincides with the joint distribution of
%
%e1.6 ###
\begin{equation}\label{multioverlap}
R^\infty_{l_1,\ldots, l_n}
=
\mathbb{E}_v  \bar{\sigma}(w,u_{l_1},v)\cdots\bar{\sigma}(w,u_{l_n},v)
\end{equation}
for $\bar{\sigma} = \bar{\sigma}_\mu$, for all $n\geq1$ and all
$l_1,\ldots, l_n\geq1$,
where $\mathbb{E}_v$ is the expectation in the last coordinate $v$
only. For $n=2$, the corresponding
quantity
%
%e1.7 ###
\begin{equation}\label{overlaps}
R^\infty_{l,l'}=\mathbb{E}_v  \bar{\sigma}(w,u_l,v) \bar{\sigma
}(w,u_{l'},v)
\end{equation}
is the\vspace*{1pt} asymptotic version of the overlap $N^{-1}\sum_{i\leq N}\sigma
_i^l\sigma_i^{l'}$.
With these notations it is clear that the famous Parisi ultrametricity
conjecture, which says
that $R_{2,3}^\infty\geq\min(R_{1,2}^\infty, R_{1,3}^\infty)$ with
probability one,
can be expressed in terms of $\bar{\sigma}_\mu$ by saying that for
all $w\in[0,1]$
the family of functions $v\to\bar{\sigma}_\mu(w,u,v)$ parametrized
by $u\in[0,1]$
is ultrametric in $\LL^2([0,1],dv)$.

An ultimate goal would be to show that the set of possible limits $\mu
\in\M$ and their
representations $\sigma_\mu$ are described by the Parisi ultrametric ansatz.
Even though this goal is out of reach at the moment, in the setting of
the Sherrington--Kirkpatrick
and diluted models we will obtain several results which demonstrate
that the point of view
based on the Aldous--Hoover representation (\ref{sigma}) provides a
useful framework
for studying the asymptotic behavior of these models.
First, we will narrow down possible limits $\M$ to some well-defined
class of distributions
$\M_{\mathrm{inv}}$ that will be described via invariance and
self-consistency
equations on~$\sigma_\mu$.
The proof of these invariance properties will be based on some standard
cavity computations; however, justification of these computations will
rely on certain properties of convergence
of measures $\mu_N$ that are not intuitive or, at least, do not easily
follow from known results.
In both types of models we will introduce a novel perturbation of
the Hamiltonian that will force the sequence $(\mu_N)$ to satisfy
these properties, and the ideas behind these perturbations will
constitute the main technical contribution of the paper.

Besides giving some constructive description of possible limits
$\mathcal M$, the invariance equations
will play a significant role in other ways. First, using these
equations we will be able
to prove representations for the limit of the free energy $F_N=N^{-1}
\mathbb{E}\log Z_N$ in terms
of $\sigma_\mu$ for $\mu\in\M_{\mathrm{inv}}$ which will automatically
coincide with the corresponding
Parisi formulas for the free energy if one can show that all measures
in $\mathcal{M}_{\mathrm{inv}}$ satisfy
the predictions of the Parisi ansatz.
These representations, proved in Sections~\ref{SecLB},~\ref{SecUB}
for diluted models and
in Sections~\ref{SecLBSK},~\ref{SecUBSK} for the
Sherrington--Kirkpatrick model, will arise
from an application of the Aizenman--Sims--Starr scheme introduced in~\cite{AS2} and,
what is crucial, thanks to the invariance equations we will only use
this scheme with
one cavity coordinate whereas all previous applications of this scheme (e.g.,
in~\cite{AS2,DeSanctis} or~\cite{PCS}) only worked when the
number of cavity
coordinates goes to infinity.

In the setting of the Sherrington--Kirkpatrick model we will utilize a
Gaussian nature of
the Hamiltonian to give other important applications of the invariance
properties of
$\mu\in\mathcal{M}$. First, we will prove in Theorem~\ref{ThSKo}
below that the joint distributions
of all spins, and thus measure $\mu$, are completely determined by the
joint distribution of
the overlaps (\ref{overlaps}). Then in Section~\ref{SecParisiA} we
will show that all limits
$\mu\in\mathcal M$ that satisfy the Parisi ultrametricity hypothesis
correspond to $\sigma_\mu$
given by certain specific realizations of the Ruelle probability
cascades. This means that,
under ultrametricity, we obtain a more detailed asymptotic description
of the model which
includes the joint distribution of all spins or multi-overlaps and not
only overlaps,
as in the usual description of the Parisi ansatz.
Motivated by this special form of ${\sigma}_\mu$ in the
Sherrington--Kirkpatrick model,
in the second part of Section~\ref{SecParisiA} we will try to
formulate a more general Parisi ansatz
expected to hold in the diluted models in terms of the Aldous--Hoover
representation (\ref{sigma}).

Finally, we would like to mention recent work~\cite{AC} where the
authors study asymptotic
behavior of spin glass models in the framework of random overlap
structures, or ROSts,
which in our notation correspond to the $\LL^2([0,1],dv)$ structure
of the family
of functions $v\to\bar{\sigma}_\mu(w,u,v)$. They obtain a number of
interesting properties
of ROSts and prove several results which are similar in spirit to ours,
for example, the Parisi formula
in the Sherrington--Kirkpatrick model under the assumption of ultrametricity.

%s1.1 ###
\subsection{Diluted models}
To illustrate the main new ideas we will start with the case of the
diluted models
where many technical details will be simpler.
We will consider the following class of diluted models as in~\cite{PT}.
Let $p\geq2$ be an \textit{even} integer, and let $\alpha>0$. Consider a
random function $\theta\dvtx \{-1,+1\}^p \to\Reals$
and a sequence $(\theta_k)_{k\geq1}$ of independent copies of $\theta$.
Consider an i.i.d. sequence of indices $(i_{l,k})_{l,k\geq1}$ with
uniform distribution
on $\{1,\ldots,N\}$, and let $\pi(\alpha N)$ be a Poisson r.v. with
mean $ \alpha N$.
Let us define the Hamiltonian $H_N(\vsi)$ on
$\Sigma_N$ by
%
%e1.8 ###
\begin{equation}\label{Ham}
-H_N(\vsi)=\sum_{k\leq\pi(\alpha N)}\theta_k(\sigma_{i_{1,k}},
\ldots,\sigma_{i_{p,k}}).
\end{equation}
Clearly, any such model has symmetry between sites.
We will make the following assumptions on the random function $\theta$.
We assume that there exists a random function $f\dvtx\{-1,+1\}\to\Reals$
[i.e., $f(\sigma) = f' + f''\sigma$ for some random $(f',f'')$]
such that
%
%e1.9 ###
\begin{equation}\label{condition1}
\exp\theta(\sigma_1,\ldots,\sigma_p) = a\bigl(1+ b f_1(\sigma_1)\cdots
f_p(\sigma_p)\bigr),
\end{equation}
where $f_1,\ldots,f_p$ are independent copies of $f$,
$b$ is a r.v. independent of $f_1,\ldots,f_p$ that satisfies the condition
%
%e1.10 ###
\begin{equation}\label{condition2}
\forall n\geq1 \qquad\mathbb{E}(-b)^n \geq0,
\end{equation}
and $a$ is an arbitrary r.v. such that $\mathbb{E}|{\log a}|<\infty$.
Finally, we assume that
%
%e1.11 ###
\begin{equation}\label{condition3}
|bf_1(\sigma_1)\cdots f_p(\sigma_p)|< 1 \qquad\mbox{ a.s.},
\end{equation}
and $\theta$ satisfies some mild integrability conditions
%
%e1.12 ###
\begin{equation}\label{condition4}
-\infty< \mathbb{E}\min_{\sigma} \theta(\sigma_1,\ldots,\sigma
_p),\qquad
\mathbb{E}\max_{\sigma} \theta(\sigma_1,\ldots,\sigma_p)
<+\infty.
\end{equation}
Two well-known models in this class of models are the $p$-spin and
$K$-sat models.
\begin{Example}[($p$-spin model)]\label{Example1}
Consider $\beta>0$ and
a symmetric r.v. $J$. The $p$-spin model corresponds to
the choice of
\[
\theta(\sigma_1,\ldots,\sigma_p)= \beta J \sigma_1\cdots\sigma_p.
\]
Equation (\ref{condition1}) holds with
$a=\ch(\beta J)$, $b = \th(\beta J)$ and $f(\sigma)=\sigma$
and condition (\ref{condition2}) holds
since we assume that the distribution of $J$ is symmetric.
Equation (\ref{condition4}) holds if $\mathbb{E}|J|<\infty$.
\end{Example}
\begin{Example}[($K$-sat model)]\label{Example2}
Consider $\beta>0$ and a sequence of i.i.d. Bernoulli r.v.
$(J_l)_{l\geq1}$ with $\mathbb{P}(J_l=\pm1)=1/2$. The $K$-sat model
(with $K=p$) corresponds to
\[
\theta(\sigma_1,\ldots,\sigma_p)=-\beta\prod_{l\leq p} \frac
{1+J_l \sigma_l}{2}.
\]
Equation (\ref{condition1}) holds with
$a=1$, $b=e^{-\beta}-1$ and $f_l(\sigma_l)=(1+J_l\sigma_l)/2$,
and (\ref{condition2}) holds since $b<0$.
\end{Example}

It is well known that under the above conditions the sequence $NF_N$ is
super-additive,
and, therefore, the limit of $F_N$ exists; see, for example, \cite
{DeSanctis}. If we knew that $(\mu_N)$
has a unique limit, that is, $\M= \{\mu\}$, then computing the limit of
the free energy in terms of
$\sigma_\mu$ in (\ref{sigma}) would be rather straightforward as
will become clear in Section
\ref{SecLB}. However, since we do not know how to prove that $(\mu
_N)$ converges, this will create
some obstacles. Moreover, if $(\mu_{N_k})$ converges to $\mu$ over
some subsequence $(N_k)$ we do not
know how to show that $(\mu_{N_k+n})$ converges to the same limit for
a fixed shift $n\geq1$, even
though we can show that it does converge simply by treating $n$ of the
coordinates as cavity coordinates.
Even if we knew that $\mu_N$ converges, we would still like to have
some description of what the limit looks like.
To overcome some of these obstacles, we will utilize the idea of adding
a ``small'' perturbation to the Hamiltonian (\ref{Ham}) that
will not affect the limit of the free energy but at the same time
ensure that $(\mu_{N_k+n})$
and $(\mu_{N_k})$ converge to the same limit. In some sense, this is
similar to the idea of adding
$p$-spin perturbation terms in the Sherrington--Kirkpatrick model to
force the overlap distribution to satisfy the Ghirlanda--Guerra
identities~\cite{GG}; see also~\cite{DF}. The perturbation for
diluted models will be defined as follows.

Consider a sequence $(c_N)$ such that $c_N \to\infty$, $c_N/N\to0$
and $|c_{N+1}-c_N|\to0$. Consider an i.i.d. sequence of indices
$(i_{j,k,l})_{j,k,l\geq1}$ with
uniform distribution on $\{1,\ldots,N\}$, let $\pi(c_N)$ be a Poisson
r.v. with mean $c_N$,
$(\pi_l(\alpha p))$ be i.i.d. Poisson with mean $\alpha p$ and
$(\theta_{k,l})$ be a sequence
of i.i.d. copies of~$\theta$. All these random variables are assumed
to be independent
of each other and of everything else. Whenever we introduce a new
random variable,
by default it is assumed to be independent of all other random variables.
Let us define the perturbation Hamiltonian $H_N^p(\vsi)$ on
$\Sigma_N$ by
%
%e1.13 ###
\begin{equation}\label{Hampert}\quad
-H_N^p(\vsi)=\sum_{l\leq\pi(c_N)}
\log\Av_\eps
\exp\sum_{k\leq\pi_l(\alpha p)}
\theta_{k,l}(\eps, \sigma_{i_{1,k,l}}, \ldots,\sigma_{i_{p-1,k,l}}),
\end{equation}
where $\Av_\eps$ will denote uniform average over $\eps\in\{-1,+1\}$
as well as replicas $(\eps_l)$ below. Let us redefine the Hamiltonian
in (\ref{Ham}) by
%
%e1.14 ###
\begin{equation}\label{Hampluspert}
-H_N(\vsi)=\sum_{k\leq\pi(\alpha N)}\theta_k(\sigma_{i_{1,k}},
\ldots,\sigma_{i_{p,k}}) - H_{N}^p(\vsi),
\end{equation}
and from now on we assume that $(\mu_N)$ and $\M$ are defined for
this perturbed Hamiltonian.
Obviously, condition (\ref{condition4}) implies that the perturbation
term does not
affect the limit of free energy since $c_N=o(N)$. The benefits of
adding this perturbation term
will first appear in Lemma~\ref{LemShift2} below where it will be
shown that thanks to this term
$(\mu_{N_k})$ and $(\mu_{N_k+n})$ converge to the same limit for any
fixed shift $n\geq1$.
Another important consequence will appear in Theorem~\ref{ThSC} below
where the perturbation
will force the limiting distributions $\mu\in\M$ to satisfy some
important invariance properties
that will play crucial role in the proof of the representation for the
free energy in Theorem~\ref{ThFE}.

Let us introduce some notations. We will usually work with $\sigma_\mu
$ for a fixed distribution
$\mu\in\M$ so for simplicity of notation we will omit subscript $\mu
$ and simply write $\sigma$.
Let $(v_{i_1,\ldots, i_n}), (x_{i_1,\ldots,i_n})$ be i.i.d. sequences
uniform on $[0,1]$
for $n\geq1$ and $i_1,\ldots, i_n\geq1$, and let
%
%e1.15 ###
\begin{equation}\label{s}
s_{i_1,\ldots, i_n} = \sigma(w,u, v_{i_1,\ldots, i_n},x_{i_1,\ldots
, i_n}).
\end{equation}
The role of multi-indices $(i_1,\ldots,i_n)$ will be simply to select
various subsets of array
(\ref{sigma}) with disjoint coordinate indices $i$ without worrying
about how to enumerate them.
Let $(\theta_{i_1,\ldots, i_n})$ be the copies of random function
$\theta$ independent over
different sets of indices. In addition, let $\hat{v}, \hat{x}, \hat
{\theta}$ be independent copies
of the above sequences, and let
%
%e1.16 ###
\begin{equation}\label{hats}
\hat{s}_{i_1,\ldots, i_n} = \sigma(w,u, \hat{v}_{i_1,\ldots, i_n},
\hat{x}_{i_1,\ldots, i_n}).
\end{equation}
Notice that we keep the same $w$ and $u$ in both $s$ and $\hat{s}$.
Throughout the paper let us denote by $\pi(\lambda)$ Poisson random
variables with mean
$\lambda$ which will always be independent from all other random
variables and from each other.
For example, if we write $\pi(\alpha)$ and $\pi(\beta)$, we assume
them to be independent
even if $\alpha=\beta$. Let $(\pi_j(\lambda))$ be independent
copies of these r.v. for $j\geq1$.
Let
%
%e1.17 ###
\begin{equation}\label{Ai}
A_i(\eps)=\sum_{k\leq\pi_i(p\alpha)} \theta_{k,i}(\eps,
s_{1,i,k},\ldots,s_{p-1,i,k})
\end{equation}
for $i\geq1$ and $\eps\in\{-1,+1\}$, and let
%
%e1.18 ###
\begin{equation}\label{Bi}
B_i = \sum_{k\leq\pi_i((p-1)\alpha)} \hat{\theta}_{k,i}(\hat
{s}_{1,i,k},\ldots,\hat{s}_{p,i,k}).
\end{equation}
We will express invariance and self-consistency properties of
distributions $\mu\in\M$
in terms of equations for the joint moments of arbitrary subset of spins
in the array~(\ref{sigma}). Take arbitrary $n, m, q, r\geq1$ such
that $n\leq m$.
In the equations below, the index $q$ will correspond to the number of
replicas selected,
$m$ will be the total number of coordinates and $n$ the number of
cavity coordinates considered
and $r$ will be the number of perturbation terms of certain type.
For each replica index $l\leq q$ we consider an arbitrary subset of coordinates
$C_l\subseteq\{1,\ldots, m\}$ and split them into the cavity and
noncavity coordinates
%
%e1.19 ###
\begin{equation}\label{C12}
C_l^1 = C_l\cap\{1,\ldots, n\},\qquad
C_l^2=C_l\cap\{n+1,\ldots,m\}.
\end{equation}
Let $\mathbb{E}'$ denote the expectation in $u$ and in sequences $x$
and $\hat{x}$, and let
%
%e1.20 ###
\begin{equation}\label{Ul}
U_l = \mathbb{E}' \Av_\eps\prod_{i\in C_l^1} \eps_i \exp\sum
_{i\leq n} A_i(\eps_i)
\prod_{i\in C_l^2} s_i \exp\sum_{k\leq r} \hat{\theta}_k(\hat
{s}_{1,k},\ldots,\hat{s}_{p,k})
\end{equation}
and
%
%e1.21 ###
\begin{equation}\label{Vl}
V = \mathbb{E}' \Av_\eps\exp\sum_{i\leq n} A_i(\eps_i)
\exp\sum_{k\leq r} \hat{\theta}_k(\hat{s}_{1,k},\ldots,\hat{s}_{p,k}).
\end{equation}
Then the following holds.
\begin{theorem}\label{ThSC}
For any limiting distribution $\mu\in\M$ and $\sigma= \sigma_\mu
$, we have
%
%e1.22 ###
\begin{equation}\label{SC}
\mathbb{E}\prod_{l\leq q} \prod_{i\in C_l}s_i^l
=\mathbb{E}\prod_{l\leq q} \mathbb{E}' \prod_{i\in C_l}s_i
=\mathbb{E}\frac{\prod_{l\leq q}U_l}{V^q}.\vadjust{\goodbreak}
\end{equation}
\end{theorem}

We will say a few words about various interpretations of (\ref{SC})
below, but first let us describe
the promised representation for the free energy.
Let $\M_{\mathrm{inv}}$ denote the set of distributions of exchangeable arrays
generated by functions
$\sigma\dvtx[0,1]^4\to\{-1,+1\}$ as in (\ref{sigma}) that satisfy
invariance equations (\ref{SC}) for all
possible choices of parameters. Theorem~\ref{ThSC} proves that $\M
\subseteq\M_{\mathrm{inv}}$. Let
\[
A(\eps)
=
\sum_{k\leq\pi(p \alpha)} \theta_{k}(\eps, s_{1,k},\ldots,s_{p-1,k})
\]
for $\eps\in\{-1,+1\}$,
\[
B
=
\sum_{k\leq\pi((p-1) \alpha)} \theta_{k}(s_{1,k},\ldots,s_{p,k})
\]
and let
%
%e1.23 ###
\begin{equation}\label{CalP}
\P(\mu)=
\log2 + \mathbb{E}\log\mathbb{E}' \Av_{\eps}
\exp A(\eps)
-\mathbb{E}\log\mathbb{E}' \exp B.
\end{equation}
The following representation holds.
\begin{theorem}\label{ThFE} We have
%
%e1.24 ###
\begin{equation}\label{FE}
\lim_{N\to\infty} F_N =
\inf_{\mu\in\M} \P(\mu) =
\inf_{\mu\in\M_{\mathrm{inv}}} \P(\mu).
\end{equation}
\end{theorem}

One can simplify the last term in (\ref{CalP}) since we will show at
the end of Section~\ref{SecUB} that
%
%e1.25 ###
\begin{equation}\label{Plast}
\mathbb{E}\log\mathbb{E}' \exp B
=(p-1)\alpha\mathbb{E}\log\mathbb{E}' \exp\theta(s_1,\ldots,s_p)
\end{equation}
for $\mu\in\M_{\mathrm{inv}}$. To better understand (\ref{SC}) let us
describe several special cases.
Let us define
%
%e1.26 ###
\begin{equation}\label{Aiav}
A_i = \log\Av_\eps\exp A_i(\eps).
\end{equation}
First, if we set $r=0$ and let sets $C_l$ be such that $C_l \subseteq\{
n+1,\ldots, m\}$ for all
$l\leq q$, then (\ref{SC}) becomes
%
%e1.27 ###
\begin{equation}\label{ASC}
\mathbb{E}\prod_{l\leq q}\mathbb{E}' \prod_{i\in C_l} s_i
=
\mathbb{E}\frac{\prod_{l\leq q} \mathbb{E}' \prod_{i\in C_l} s_i
\exp\sum_{i\leq n} A_i}
{(\mathbb{E}' \exp\sum_{i\leq n} A_i)^q}.
\end{equation}
On the other hand, if we set $n=0$, then (\ref{SC}) becomes
%
%e1.28 ###
\begin{equation}\label{preBSC}
\mathbb{E}\prod_{l\leq q}\mathbb{E}'\prod_{i\in C_l} s_i
=
\mathbb{E}\frac{\prod_{l\leq q} \mathbb{E}' \prod_{i\in C_l} s_i
\exp\sum_{i\leq r} \hat{\theta}_i(\hat{s}_{1,i},\ldots,\hat{s}_{p,i})}
{(\mathbb{E}' \exp\sum_{i\leq r} \hat{\theta}_i(\hat
{s}_{1,i},\ldots,\hat{s}_{p,i}))^q}.
\end{equation}
These equations can be interpreted as the invariance of the
distribution of $(s_i^l)$
under various changes of density, and they will both play an important role
in the proof of Theorem~\ref{ThFE}. Another consequence of (\ref{SC})
are the following
self-consistency equations for the distribution of\vadjust{\goodbreak} spins. Let us set
$r=0$ and $n=m$.
Let
\[
s_i^A = \frac{\Av_\eps  \eps\exp A_i(\eps)}{\Av_\eps\exp
A_i(\eps)}.
\]
Then (\ref{SC}) becomes
%
%e1.29 ###
\begin{equation}\label{ASCSC}
\mathbb{E}\prod_{l\leq q}\mathbb{E}' \prod_{i\in C_l} s_i
=
\mathbb{E}\frac{\prod_{l\leq q} \mathbb{E}' \prod_{i\in C_l} s_i^A
\exp\sum_{i\leq n} A_i}
{ (\mathbb{E}' \exp\sum_{i\leq n} A_i)^q}.
\end{equation}
This means that the distribution of spins $(s_i^l)$ coincides with the
distribution of ``new'' spins
$(s_i^{A,l})$ under a certain change of density. Even though we cannot
say more about the role
(\ref{ASCSC}) might play in the diluted models, its analog in the
Sherrington--Kirkpatrick model
will play a very important role in proving that the joint overlap
distribution under $\mu$ determines
$\mu$ and in constructing the explicit formula for $\bar{\sigma}$
under the Parisi ultrametricity
hypothesis.

It will become clear from the arguments below that, in essence, the
representation (\ref{FE}) is
the analog of the Aizenman--Sims--Starr scheme in the
Sherrington--Kirkpatrick model~\cite{AS2}
with one cavity coordinate. Previous applications of this scheme (e.g.,
in~\cite{AS2,DeSanctis} or~\cite{PCS}) only worked when the number of cavity
coordinates goes to infinity,
since considering one cavity coordinate in general yields only a lower
bound on the free energy.
This lower bound expressed in terms of the generic functional order
parameter $\sigma_\mu$ will
be proved in Section~\ref{SecLB}. Then the main new ideas of the
paper---the roles played by
the perturbation Hamiltonian (\ref{Hampert}) and the consequent
invariance in (\ref{SC})---will help us justify that this lower bound
is exact and, moreover, represent it via a well-defined
family $\M_{\mathrm{inv}}$. First, following the arguments in~\cite{FL,PT},
in Section~\ref{SecUB}
we will prove a corresponding Franz--Leone type upper bound which will
depend on an arbitrary function $\sigma$ that defines an exchangeable
array as in (\ref{sigma}). For a general $\sigma$,
this upper bound will depend on $N$. However, we will show that for
$\sigma_\mu$ for $\mu\in
\M_{\mathrm{inv}}$ the invariance of Theorem~\ref{ThSC} implies that the upper
bound is independent
of $N$ and matches the lower bound. This is the main point where the
invariance properties
will come into play. The same ideas will work in the
Sherrington--Kirkpatrick model with the
appropriate choice of the perturbation Hamiltonian.

%s1.2 ###
\subsection{The Sherrington--Kirkpatrick model}
Let us consider mixed $p$-spin Sherrington--Kirkpatrick Hamiltonian
%
%e1.30 ###
\begin{equation}\label{HamSK}
-H_N(\vsi)
=- \sum_{p\geq1} \beta_p H_{N,p}(\vsi),
\end{equation}
where
%
%e1.31 ###
\begin{equation}\label{HamSKp}
- H_{N,p}(\vsi)
=
\frac{1}{N^{(p-1)/2}}
\sum_{1\leq i_1,\ldots,i_p\leq N}g_{i_1,\ldots,i_p} \sigma
_{i_1}\cdots\sigma_{i_p},\vadjust{\goodbreak}
\end{equation}
the sum is over $p=1$ and even $p\geq2$ and
$(g_{i_1,\ldots,i_p})$ are standard Gaussian independent for all
$p\geq1$ and
all $(i_1,\ldots,i_p)$. The covariance of (\ref{HamSK}) is given
by\looseness=-1
%
%e1.32 ###
\begin{equation}
\mathbb{E}H_N(\vsi^1) H_N(\vsi^2) = N\xi(R_{1,2}),
\end{equation}\looseness=0
where $\xi(x)=\sum_{p\geq1}\beta_p^2 x^p$, and we assume that the
sequence $(\beta_p)$ satisfies $\sum_{p\geq1} 2^p \beta_p^2<\infty$.
Let us\vspace*{1pt} start by introducing the analog of the perturbation
Hamiltonian (\ref{Hampert}) for the Sherrington--Kirkpatrick model.
Consider independent Gaussian processes $G_{\xi'}(\vsi)$ and
$G_{\theta}(\vsi)$ on $\Sigma_N=\{-1,+1\}^N$ with covariances
%
%e1.33 ###
\begin{equation}\label{covSK}\quad
\mathbb{E}G_{\xi'}(\vsi^1) G_{\xi'}(\vsi^2)
=\xi'(R_{1,2}),\qquad
\mathbb{E}G_{\theta}(\vsi^1) G_{\theta}(\vsi^2)
=\theta(R_{1,2}),
\end{equation}
where $\theta(x)=x\xi'(x)-\xi(x)$, and let $G_{\xi',k}(\vsi)$ and
$G_{\theta,k}(\vsi)$ be
their independent copies for $k\geq1$. For $(c_{N})$ as above, let us
add the following
perturbation to the Hamiltonian (\ref{HamSK}):
%
%e1.34 ###
\begin{equation}\label{HampertSK}
-H_N^{p}(\vsi)
=\sum_{k\leq\pi(c_N)}\log\ch  G_{\xi',k}(\vsi)
+\sum_{k\leq\pi'(c_N)} G_{\theta,k}(\vsi),
\end{equation}
where $\pi(c_N)$ and $\pi'(c_N)$ are independent Poisson random
variables with means~$c_N$.
Clearly, this Hamiltonian does not affect the limit of the free energy
since $c_N=o(N)$.
We will see that this choice of perturbation ensures the same nice
properties of
convergence as the perturbation (\ref{Hampert}) in the setting of the
diluted models.
As a consequence, we will get the following analog of the invariance of
Theorem~\ref{ThSC}.
Given a measurable function $\bar{\sigma}\dvtx[0,1]^3\to[-1,1]$, for
any $w\in[0,1]$, let
$g_{\xi'}(\bar{\sigma}(w,u,\cdot))$ be a Gaussian process indexed
by functions $v\to\bar{\sigma}(w,u,\cdot)$
for $u\in[0,1]$ with covariance
%
%e1.35 ###
\begin{equation}\label{gxicov}\qquad
\mbox{Cov} ( g_{\xi'}(\bar{\sigma}(w,u,\cdot)), g_{\xi
'}(\bar{\sigma}(w,u',\cdot)) )
=
\xi'(\mathbb{E}_v \bar{\sigma}(w,u,v) \bar{\sigma
}(w,u',v))
\end{equation}
and $g_{\theta}(\bar{\sigma}(w,u,\cdot))$ be a Gaussian process
independent of $g_{\xi'}(\bar{\sigma}(w,u,\cdot))$
with covariance
%
%e1.36 ###
\begin{equation}\label{gthetacov}\qquad
\mbox{Cov} ( g_{\theta}(\bar{\sigma}(w,u,\cdot)), g_{\theta
}(\bar{\sigma}(w,u',\cdot)) )
=\theta(\mathbb{E}_v \bar{\sigma}(w,u,v) \bar{\sigma}
(w,u',v)).
\end{equation}
Let us consider independent standard Gaussian random variables $z$ and
$z'$ and define
%
%e1.37 ###
\begin{equation}\label{Gxi}\qquad
G_{\xi'}(\bar{\sigma}(w,u,\cdot)) = g_{\xi'}(\bar{\sigma
}(w,u,\cdot))
+z \bigl(\xi'(1)- \xi'(\mathbb{E}_v \bar{\sigma}(w,u,v)^2)\bigr)^{1/2}
\end{equation}
and
%
%e1.38 ###
\begin{equation}\label{Gtheta}\qquad
G_{\theta}(\bar{\sigma}(w,u,\cdot)) = g_{\theta}(\bar{\sigma
}(w,u,\cdot))
+ z' \bigl(\theta(1)- \theta(\mathbb{E}_v \bar{\sigma
}(w,u,v)^2)\bigr)^{1/2}.
\end{equation}
For simplicity of notation we will keep the dependence of $G_{\xi'}$
and $G_\theta$ on $z$
or $z'$ implicit. Let $G_{\xi', i}$ and $G_{\theta,i}$ be independent
copies of these processes.
Random variables $z$, and $z'$ will play the role of replica variables
similarly to $u$ and
for this reason in the Sherrington--Kirkpatrick model we will denote by
$\mathbb{E}'$ the expectation
in $u$, $z$ and $z'$. The main\vadjust{\goodbreak} purpose of introducing the second term
in (\ref{Gxi}) and
(\ref{Gtheta}) is to match the variances of these Gaussian processes,
$\xi'(1)$ and $\theta(1)$,
to variances in (\ref{covSK}) for $\vsi^1=\vsi^2$.

As in the setting of diluted models, consider arbitrary $n, m, q, r\geq
1$ such that $n\leq m$.
For each $l\leq q$ consider an arbitrary subset $C_l\subseteq\{
1,\ldots, m\}$, and let
$C_l^1$ and $C_l^2$ be defined as in (\ref{C12}). Let $\bar{\sigma
}_i = \bar{\sigma}(w,u,v_i)$.
For $l\leq q$ define
%
%e1.39 ###
\begin{equation}
U_l = \mathbb{E}'
\prod_{i\in C_l^1}\th  G_{\xi',i}(\bar{\sigma}(w,u,\cdot)) \prod
_{i\in C_l^2} \bar{\sigma}_i
  \E_{n,r},
\end{equation}
where
%
%e1.40 ###
\begin{equation}\quad
\E_{n,r}= \exp\biggl( \sum_{i\leq n} \log\ch  G_{\xi',i}(\bar
{\sigma}(w,u,\cdot))
+ \sum_{k\leq r} G_{\theta,k}(\bar{\sigma}(w,u,\cdot))
\biggr),
\end{equation}
and let $V=\mathbb{E}' \E_{n,r}$. If $\M$ denotes the set of
possible limits of $\mu_N$
corresponding to the Hamiltonian (\ref{HamSK}) perturbed by (\ref
{HampertSK}), then the following holds.

\begin{theorem}\label{ThSCSK}
For any $\mu\in\M$ and $\bar{\sigma}=\bar{\sigma}_\mu$ we have
%
%e1.41 ###
\begin{equation}\label{SCSK}
\mathbb{E}\prod_{l\leq q} \mathbb{E}' \prod_{i\in C_l}\bar{\sigma}_i
=\mathbb{E}  \frac{\prod_{l\leq q}U_l}{V^q}.
\end{equation}
\end{theorem}

Let $\M_{\mathrm{inv}}$ be the family of distributions defined by the
invariance properties (\ref{SCSK}),
so that Theorem~\ref{ThSCSK} proves that $\M\subseteq\M_{\mathrm{inv}}$. If
we define
%
%e1.42 ###
\begin{eqnarray}\label{CalPSK}
\P(\mu)
&=&\log2+\mathbb{E}\log \mathbb{E}' \ch
G_{\xi'}(\bar{\sigma}_\mu(w,u,\cdot))\nonumber\\[-8pt]\\[-8pt]
&&{}-\mathbb{E}\log \mathbb{E}' \exp
G_{\theta}(\bar{\sigma}_\mu(w,u,\cdot)),\nonumber
\end{eqnarray}
then we have the following representation for the free energy in the
Sherrington--Kirkpatrick model.
\begin{theorem}\label{ThFESK}
We have
%
%e1.43 ###
\begin{equation}\label{Parisi}
\lim_{N\to\infty} F_N = \inf_{\mu\in\M} \P(\mu)
= \inf_{\mu\in\M_{\mathrm{inv}}}\P(\mu).
\end{equation}
\end{theorem}

As in the case of diluted models above, let us describe several special
cases of (\ref{SCSK}).
If $r=0$ and sets $C_l$ are such that $C_l \subseteq\{n+1,\ldots, m\}
$ for all
$l\leq q$, then (\ref{SCSK}) becomes
%
%e1.44 ###
\begin{equation}\label{ASCSK}
\mathbb{E}\prod_{l\leq q}\mathbb{E}' \prod_{i\in C_l} \bar{\sigma}_i
=\mathbb{E}  \frac{\prod_{l\leq q} \mathbb{E}' \prod_{i\in C_l}
\bar{\sigma}_i \prod_{i\leq n}
\ch  G_{\xi',i}(\bar{\sigma}(w,u,\cdot))}
{ (\mathbb{E}' \prod_{i\leq n} \ch  G_{\xi',i}(\bar{\sigma
}(w,u,\cdot)))^q}.
\end{equation}
If we set $n=0$, then (\ref{SCSK}) becomes
%
%e1.45 ###
\begin{equation}\label{preBSCSK}\quad
\mathbb{E}\prod_{l\leq q}\mathbb{E}'\prod_{i\in C_l} \bar{\sigma}_i
=
\mathbb{E}  \frac{\prod_{l\leq q} \mathbb{E}' \prod_{i\in C_l}
\bar{\sigma}_i
\exp\sum_{k\leq r} G_{\theta,k}(\bar{\sigma}(w,u,\cdot))}
{(\mathbb{E}' \exp\sum_{k\leq r} G_{\theta,k}(\bar{\sigma
}(w,u,\cdot)))^q}.\vadjust{\goodbreak}
\end{equation}
Again, these equations can be interpreted as the invariance of the spin
distributions under various
random changes of density. Finally, if we set $r=0$ and $n=m$, then
(\ref{SCSK}) becomes
%
%e1.46 ###
\begin{eqnarray}\label{Invar}\qquad
&&\mathbb{E}\prod_{l\leq q}\mathbb{E}' \prod_{i\in C_l}
\bar{\sigma}_i\nonumber\\[-8pt]\\[-8pt]
&&\qquad=
\mathbb{E}  \frac{\prod_{l\leq q} \mathbb{E}' \prod_{i\in C_l}
\th  G_{\xi',i}(\bar{\sigma}(w,u,\cdot))
\prod_{i\leq n} \ch  G_{\xi',i}(\bar{\sigma}(w,u,\cdot))}
{ (\mathbb{E}' \prod_{i\leq n} \ch  G_{\xi',i}(\bar{\sigma
}(w,u,\cdot)))^q}.\nonumber
\end{eqnarray}
The meaning of this self-consistency equation is that the joint
distribution of spins generated
by a function $\bar{\sigma}(w,u,v)$ coincides with the distribution
of spins generated by
$\th  G_{\xi'}(\bar{\sigma}(w,u,\cdot))$ under a properly
interpreted random change of density,
and we will discuss this interpretation in more detail below under the
Parisi ultrametricity hypothesis.
The choice of parameters in (\ref{Invar}), most importantly $n=m$,
will be the key to the following special property of the
Sherrington--Kirkpatrick model.
\begin{theorem}\label{ThSKo}
For any $\mu\in\M_{\mathrm{inv}}$, the joint distribution of
$(R_{l,l'}^\infty)_{l,l'\geq1}$ defined
in (\ref{overlaps}) for $\bar{\sigma}= \bar{\sigma}_\mu$ uniquely
determines $\mu$
and thus the joint distribution of all multi-overlaps.
\end{theorem}

The fact that the joint distribution of overlaps determines $\mu$
leads to a natural addition to
the statement of Theorem~\ref{ThFESK}. It will be clear early in the
proof of Theorem~\ref{ThFESK}
that $\P(\mu)$ for $\mu\in\M$ depends only on the distribution of
the array (\ref{overlaps})
for $\bar{\sigma}= \bar{\sigma}_\mu$, and, as a result, one can
express the free energy in
(\ref{Parisi}) as the infimum over a family of measures $\M_{\mathrm{inv}}'$
defined completely in terms
of the invariance of the joint overlap distribution and such that $\M
_{\mathrm{inv}}\subseteq\M_{\mathrm{inv}}'$.
For this purpose one does not need the self-consistency
part of the equations (\ref{SCSK}), so we will only use the case when
$C_l^2=C_l$ in (\ref{C12})
for all $l$. Let us consider processes $G_{\xi'}$ and $G_{\theta}$ in
(\ref{Gxi}), (\ref{Gtheta}) defined
in terms of replicas $(u_l)$, $(z_l)$ and $(z_l')$ of $u$, $z$ and
$z'$, namely,
%
%e1.47 ###
\begin{equation}\label{Gxil}\
G_{\xi'}(\bar{\sigma}(w,u_l,\cdot)) = g_{\xi'}(\bar{\sigma
}(w,u_l,\cdot))
+ z_l \bigl(\xi'(1)- \xi'(\mathbb{E}_v \bar{\sigma
}(w,u_l,v)^2)\bigr)^{1/2}\hspace*{-28pt}
\end{equation}
and
%
%e1.48 ###
\begin{equation}\label{Gthetal}
G_{\theta}(\bar{\sigma}(w,u_l,\cdot)) = g_{\theta}(\bar{\sigma
}(w,u_l,\cdot))
+ z_l' \bigl(\theta(1)- \theta(\mathbb{E}_v \bar{\sigma
}(w,u_l,v)^2)\bigr)^{1/2}.\hspace*{-28pt}
\end{equation}
Let $F=F((R^\infty_{l,l'})_{l,l'\leq q})$ be an arbitrary continuous
function of the overlaps on
$q$ replicas. Let
%
%e1.49 ###
\begin{equation}
U = \mathbb{E}' F \prod_{l\leq q} \exp\biggl( \sum_{i\leq n} \log
\ch  G_{\xi',i}(\bar{\sigma}(w,u_l,\cdot))
+\sum_{k\leq r} G_{\theta,k}(\bar{\sigma}(w,u_l,\cdot))
\biggr).\hspace*{-35pt}
\end{equation}
Then the condition
%
%e1.50 ###
\begin{equation}\label{SCSKR}
\mathbb{E}F = \mathbb{E}(U/V^q)\vadjust{\goodbreak}
\end{equation}
for all $q,n,r$ and all continuous bounded functions $F$ defines the
family $\M_{\mathrm{inv}}'$. Equation (\ref{SCSKR}) is obviously implied by
(\ref{SCSK}) which contains the case of polynomial $F$ simply by
making sure that $C_l^2=C_l$, so $\M_{\mathrm{inv}}\subseteq\M_{\mathrm{inv}}'$. Then
one can add
%
%e1.51 ###
\begin{equation}\label{Parisimore}
\lim_{N\to\infty} F_N = \inf_{\mu\in\M_{\mathrm{inv}}'} \P(\mu)
\end{equation}
to the statement of Theorem~\ref{ThFESK}. This together with Theorem
\ref{ThSKo} shows that
in the Sherrington--Kirkpatrick model the role of the order parameter
is played by the joint distribution
of overlaps rather than the joint distribution of all multi-overlaps or
the generic functional order
parameter $\bar{\sigma}_\mu$. This gives an idea about how close
this point of view takes us to
the Parisi ansatz~\cite{Parisi} where the order parameter is the
distribution of one overlap.
Since we can always ensure that the Ghirlanda--Guerra identities \cite
{GG} hold by adding
a mixed $p$-spin perturbation term [see (\ref{HampertmixSK}) below],
the remaining gap is
the ultrametricity of the overlaps, since it is well known that the
Ghirlanda--Guerra identities
and ultrametricity determine the joint distribution of overlaps from
the distribution of one overlap; see, for example,~\cite{BR} or \cite
{Bovier}. If one can generalize the results in~\cite{PGG} and
\cite{Tal-New} to show that the Ghirlanda--Guerra identities always
imply ultrametricity,
(\ref{Parisi}) would coincide with the Parisi formula proved in~\cite{T-P}.

\textit{The Ghirlanda--Guerra identities and stochastic stability}.
Let us mention that the Ghirlanda--Guerra identities and stochastic stability
can also be expressed in terms of the generic functional order
parameter $\bar{\sigma}$.
We will use a version of both properties in the formulation proved in
\cite{Tal-New}.
Let us now consider a different perturbation term
%
%e1.52 ###
\begin{equation}\label{HampertmixSK}
H_N^\delta(\vsi)
=\delta_N \sum_{p\geq1} \beta_{N,p} H_{N,p}'(\vsi),
\end{equation}
where
%
%e1.53 ###
\begin{equation}\label{HampertmixSKp}
-H_{N,p}'(\vsi)
=\frac{1}{N^{(p-1)/2}}
\sum_{1\leq i_1,\ldots,i_p\leq N}g_{i_1,\ldots,i_p}' \sigma
_{i_1}\cdots\sigma_{i_p}
\end{equation}
are independent copies of (\ref{HamSKp}). When $\delta_N\to0$ this
perturbation
term is of smaller order than (\ref{HamSK}) and does not affect the
limit of the free energy.
However, the arguments in the proof of the Ghirlanda--Guerra identities
and stochastic stability in
\cite{Tal-New} require that $\delta_N$ does not go to zero too fast;
for example, the choice of
$\delta_N = N^{-1/16}$ works. Then, Theorem 2.5 in~\cite{Tal-New}
states that one can choose
a sequence $\bolds{\beta}_N=(\beta_{N,p})$ such that $|\beta
_{N,p}|\leq2^{-p}$ for all $N$ and
such that the following properties hold. First of all, if $\la\cdot
\ra$ is the Gibbs average
corresponding to the sum
%
%e1.54 ###
\begin{equation}\label{Hampertpert}
-H_N'(\vsi) = -H_N(\vsi) - H_N^\delta(\vsi)
\end{equation}
of the Hamiltonians (\ref{HamSK}) and (\ref{HampertmixSK}), and $F$
is a continuous function
of finitely many multi--overlaps\vadjust{\goodbreak} (\ref{multioverlapN}) on replicas
$\vsi^1,\ldots,\vsi^n$, then
the Ghirlanda--Guerra identities
%
%e1.55 ###
\begin{equation}\label{GGp}
\lim_{N\to\infty}
\Biggl|\mathbb{E}\la F R_{1,n+1}^p \ra - \frac{1}{n}  \mathbb{E}\la F\ra
  \mathbb{E}\la R_{1,2}^p \ra-
\frac{1}{n}\sum_{l=2}^{n} \mathbb{E}\la F R_{1,l}^p\ra\Biggr|=0
\end{equation}
hold for all $p\geq1$. Now, for $p\geq1$, let $G_{p}(\vsi)$ be a
Gaussian process
on $\Sigma_N$ with covariance
%
%e1.56 ###
\begin{equation}\label{covSKp}
\mathbb{E}G_{p}(\vsi^1) G_{p}(\vsi^2)
=R_{1,2}^p,
\end{equation}
and for $t>0$ let $\la\cdot\ra_t$ denote the Gibbs average
corresponding to the Hamiltonian
\[
-H_{N,t}'(\vsi) = -H_N'(\vsi) - t G_p(\vsi).
\]
Then, in addition to (\ref{GGp}), the following stochastic stability
property holds for any $t>0$:
%
%e1.57 ###
\begin{equation}\label{SS}
{\lim_{N\to\infty}}
|\mathbb{E}\la F \ra_t - \mathbb{E}\la F\ra|
=0.
\end{equation}
This property was also proved in~\cite{AC} without perturbation (\ref
{HampertmixSK}) under the condition
of differentiability of the limiting free energy.
Let $\mu_N$ be the joint distribution of spins (\ref{muN})
corresponding to the Hamiltonian
$H_N'(\vsi)$ and $\M$ be the set of all limits of $(\mu_N)$. Then
both (\ref{GGp}) and (\ref{SS})
can be expressed in the limit in terms of $\bar{\sigma} = \bar
{\sigma}_\mu$ for any $\mu\in\M$
as follows. First of all, (\ref{GGp}) becomes the exact equality in
the limit by comment above
(\ref{multioverlap}),
%
%e1.58 ###
\begin{equation}\label{GGplim}
\mathbb{E}F(R_{1,n+1}^\infty)^p = \frac{1}{n}  \mathbb{E}F
\mathbb{E}(R_{1,2}^\infty)^p
+\frac{1}{n}\sum_{l=2}^{n} \mathbb{E}F (R_{1,l}^\infty)^p.
\end{equation}
Stochastic stability (\ref{SS}) can be expressed as follows. For $w\in
[0,1]$, let
$g_{p}(\bar{\sigma}(w,u,\cdot))$ be a Gaussian process indexed by
$u\in[0,1]$ with covariance
%
%e1.59 ###
\begin{equation}\label{gxicovp}
\mbox{Cov} ( g_{p}(\bar{\sigma}(w,u,\cdot)), g_{p}(\bar
{\sigma}(w,u',\cdot)) )
=
(\mathbb{E}_v \bar{\sigma}(w,u,v) \bar{\sigma}(w,u',v))^p\hspace*{-28pt}
\end{equation}
and, as in (\ref{Gxi}), let
%
%e1.60 ###
\begin{equation}\label{Gxip}
G_{p}(\bar{\sigma}(w,u,\cdot)) = g_{p}(\bar{\sigma}(w,u,\cdot))
+z \bigl(1- (\mathbb{E}_v \bar{\sigma}(w,u,v)^2)^p \bigr)^{1/2}.
\end{equation}
Then (\ref{SS}) implies the following analog of Theorem~\ref{ThSCSK}.
\begin{theorem}\label{ThSSSK}
For any $\mu\in\M$ and $\bar{\sigma}=\bar{\sigma}_\mu$ we have
for all $p\geq1$
and $t>0$,
%
%e1.61 ###
\begin{equation}\label{SSSK}
\mathbb{E}\prod_{l\leq q} \mathbb{E}' \prod_{i\in C_l}\bar{\sigma}_i
=\mathbb{E}  \frac{\prod_{l\leq q} \mathbb{E}' \prod_{i\in
C_l}\bar{\sigma}_i \exp t G_p(\bar{\sigma}(w,u,\cdot))}
{(\mathbb{E}' \exp t G_p(\bar{\sigma}(w,u,\cdot)))^q}.
\end{equation}
\end{theorem}

The proof that (\ref{SS}) implies (\ref{SSSK}) will not be detailed
since it follows exactly the same
argument as the proof of Theorem~\ref{ThSCSK} (we will point this out
at the appropriate
step in Section~\ref{SecSCSK}).
Note that (\ref{SSSK}) is more general than (\ref{preBSCSK}), which
shows that
the invariance of Theorem~\ref{ThSCSK} is related to the stochastic
stability (\ref{SS}).\vadjust{\goodbreak}
It is interesting to note, however, that the size of the perturbation
(\ref{HampertSK}) that ensures
the invariance in (\ref{SCSK}) was of arbitrarily smaller order than
the original Hamiltonian
(\ref{HamSK}) since $c_N$ could grow arbitrarily slowly while
perturbation (\ref{HampertmixSK})
must be large enough since $\delta_N$ cannot go to zero too fast.
Moreover, the form of the perturbation (\ref{HampertSK}) plays a
crucial role in the proof of
the self-consistency part (\ref{Invar}) of equations (\ref{SCSK})
which will allow us to give
an explicit construction of the functional order parameter $\bar
{\sigma}$ below under the Parisi
ultrametricity hypothesis. The special case of the stochastic stability
(\ref{SSSK}) for the overlaps
[rather than multi-overlaps as in (\ref{SSSK})] was the starting point
of the main result
in~\cite{AA} under certain additional assumptions on $\bar{\sigma}$.

Let us make one more comment about the Ghirlanda--Guerra identities
(\ref{GGplim}) from
the point of view of the generic functional order parameter~$\bar
{\sigma}$. Equation (\ref{GGp})
always arises as a simple consequence of the following concentration
statement either for
the perturbation Hamiltonian (\ref{HampertmixSKp})
(see~\cite{Tal-New}),
%
%e1.62 ###
\begin{equation}
\lim_{N\to\infty} \mathbb{E}\biggl\la \biggl| \frac{H_{N,p}'}{N} - \mathbb{E}\biggl\la
\frac{H_{N,p}'}{N} \biggr\ra \biggr| \biggr\ra = 0
\end{equation}
or for the Hamiltonian in (\ref{HamSKp}),
%
%e1.63 ###
\begin{equation}
\lim_{N\to\infty} \mathbb{E}\biggl\la \biggl| \frac{H_{N,p}}{N} - \mathbb{E}\biggl\la
\frac{H_{N,p}}{N} \biggr\ra \biggr| \biggr\ra = 0,
\end{equation}
which was proved in~\cite{PGGmixed} for any $p$ such that $\beta_p
\not= 0$ in (\ref{HamSK})
(the case of $p=1$ was first proved in~\cite{Chatterjee}).
One can similarly encode the limiting Ghirlanda--Guerra identities
(\ref{GGplim}) as a concentration statement for the Gaussian process
$G_{p}(\bar{\sigma}(w,u,\cdot))$
in (\ref{Gxip}) as follows.
\begin{theorem}\label{GGlim}
Assuming (\ref{SSSK}), the following are equivalent:

\begin{longlist}[(2)]
\item[(1)]
the Ghirlanda--Guerra identities (\ref{GGplim}) hold;

\item[(2)] for all $p\geq1$,
%
%e1.64 ###
\begin{eqnarray}\label{GSS}
&&\mathbb{E}\frac{G_{p}(\bar{\sigma}(w,u,\cdot))^2 \exp t G_{p}(\bar
{\sigma}(w,u,\cdot))} {\mathbb{E}' \exp t
G_{p}(\bar{\sigma}(w,u,\cdot))}\nonumber\\[-8pt]\\[-8pt]
&&\qquad{} - \biggl(
\mathbb{E}\frac{G_{p}(\bar{\sigma}(w,u,\cdot)) \exp t G_{p}(\bar
{\sigma}(w,u,\cdot))} {\mathbb{E}' \exp t
G_{p}(\bar{\sigma}(w,u,\cdot))} \biggr)^2\nonumber
\end{eqnarray}
is uniformly bounded for all $t>0$, in which case it is equal to $1$.
\end{longlist}
\end{theorem}

The result will follow from a simple application of the Gaussian
integration by parts
and the main reason behind this equivalence will be very similar to the
proof of
the Ghirlanda--Guerra identities for Poisson--Dirichlet cascades in
\cite{SG2}.

%s1.3 ###
\subsection{Connections to the Parisi ansatz}\label{SecParisiA}

We will now discuss how the functional order parameter $\bar{\sigma
}(w,u,v)$ fits into the picture
of the ``generic ultrametric Parisi ansatz''\vadjust{\goodbreak} expected to hold in the
Sherrington--Kirkpatrick and diluted
models and believed to represent some kind of general principle in
other models as well. We will
begin
with the case of the Sherrington--Kirkpatrick model where the joint
distribution of the overlap array
(\ref{overlaps}) under the Parisi ultrametricity conjecture is well
understood, and we will use it to
give an explicit construction of $\bar{\sigma}(w,u,v)$. This will
serve as an illustration of a more
general case that will appear in the diluted models.

\subsubsection*{Parisi ansatz in the Sherrington--Kirkpatrick model}
Let us go back to the self-consistency equations (\ref{Invar}) and show
that they can be used
to give an explicit formula for the function $\bar{\sigma}$, or the
distribution of spins, under
the Parisi ultrametricity hypothesis and the Ghirlanda--Guerra
identities. In this section
we will assume that the reader is familiar with the Ruelle probability
cascades~\cite{Ruelle}
and refer to extensive literature on the subject for details.
Equation (\ref{overlaps}) defines some realization of the directing
measure of the overlap array in
the following sense. If we think of $\bar{\sigma}(w,u,\cdot)$ as a
function in $H = \LL^2([0,1],dv)$,
then the image of the Lebesgue measure on $[0,1]$ by the map $u\to\bar
{\sigma}(w,u,\cdot)$
defines a random probability measure $\eta_w$ on $H$. Equation (\ref
{overlaps})
states that the overlaps can be generated by scalar products in $H$ of
an i.i.d. sequence from this random measure. Any such measure $\eta_w$
defined on an arbitrary
Hilbert space is called the directing measure of the overlap array
$(R_{l,l'}^\infty)$.
It is defined uniquely up to a random isometry; see, for example, Lemma
4 in~\cite{PDS}, or in the case of discrete overlap the end of the
proof of Theorem~4 in~\cite{PGG}.
By Theorem 2 in~\cite{PGG}, the Ghirlanda--Guerra identities imply that
%
%e1.65 ###
\begin{equation}\label{qstar}
\mathbb{E}_v \bar{\sigma}(w,u,v)^2=q^*   \qquad\mbox{a.s.},
\end{equation}
where $q^*$ is the\vspace*{1pt} largest point in the support of the distribution of
$R_{1,2}^\infty$, and, therefore,
equation (\ref{Invar}) can be slightly simplified by getting rid of
the last term in~(\ref{Gxi}),
%
%e1.66 ###
\begin{eqnarray}\label{Invarg}\qquad
&&
\mathbb{E}\prod_{l\leq q}\mathbb{E}' \prod_{i\in C_l}
\bar{\sigma}_i\nonumber\\[-8pt]\\[-8pt]
&&\qquad=
\mathbb{E}  \frac{\prod_{l\leq q} \mathbb{E}' \prod_{i\in C_l}
\th  g_{\xi',i}(\bar{\sigma}(w,u,\cdot))
\prod_{i\leq n} \ch  g_{\xi',i}(\bar{\sigma}(w,u,\cdot))}
{ (\mathbb{E}' \prod_{i\leq n} \ch  g_{\xi',i}(\bar{\sigma
}(w,u,\cdot)))^q}.\nonumber
\end{eqnarray}
The key observation now is that the right-hand side of (\ref{Invarg})
does not depend on the
particular realization of the directing measure since the Gaussian
process $g_{\xi'}$ is defined
by its covariance function (\ref{gxicov}) which depends only on the
$\LL^2([0,1],dv)$ structure
of the family $\bar{\sigma}(w,u,\cdot)$.
Let us first interpret the right-hand side of (\ref{Invarg}) when the
overlap distribution is discrete,
%
%e1.67 ###
\begin{equation}\label{finite}
\mathbb{P}(R_{1,2}^\infty=q_l) = m_{l+1}-m_l\vadjust{\goodbreak}
\end{equation}
for some $0 \leq q_1< q_2 <\cdots< q_k=q^*\leq1$ and $0=m_1<\cdots
<m_k<m_{k+1}= 1$.
In this case it is well known that one directing measure of the
overlaps is given by the Ruelle
probability cascades, of course, assuming the Ghirlanda--Guerra
identities and ultrametricity
(see, e.g.,~\cite{AA,PGG,Tal-New} or~\cite{SG2})
and, therefore, $(g_{\xi',i})$ are the usual Gaussian fields
associated with the cascades. The Ruelle probability cascades is a
discrete random
measure with Poisson--Dirichlet weights $(w_\alpha)$ customarily
indexed by
$\alpha\in\Natural^k$, where $k$ is the number of atoms in (\ref{finite}),
so that the Gaussian fields are also indexed by $\alpha$, $(g_{\xi
',i}(\alpha))$.
By definition of the directing measure $\eta_w$,
the expectation $\mathbb{E}'$ in $u$ plays the role of averaging with
respect to these
weights, so that the right-hand side of (\ref{Invarg}) can be
rewritten as
%
%e1.68 ###
\begin{equation}\label{InvarDR}
\mathbb{E}\frac{\prod_{l\leq q} \sum_{\alpha} w_\alpha\prod
_{i\in C_l} \th  g_{\xi',i}(\alpha)
\prod_{i\leq n} \ch  g_{\xi',i}(\alpha)}
{ (\sum_{\alpha} w_\alpha\prod_{i\leq n} \ch  g_{\xi',i}(\alpha))^q}.
\end{equation}
This in its turn can be rewritten using well-known properties of the
Ruelle probability cascades,
in particular, Lemma 1.2 in~\cite{PTDR} which is a recursive
application of Proposition A.2
in~\cite{Bolthausen}. If we denote
\[
w_\alpha' = \frac{w_\alpha\prod_{i\leq n} \ch  g_{\xi',i}(\alpha)}
{ \sum_{\alpha} w_\alpha\prod_{i\leq n} \ch  g_{\xi',i}(\alpha)},
\]
then the point processes
%
%e1.69 ###
\begin{equation}\label{gprime}
( w_\alpha', (g_{\xi',i}(\alpha))_{i\leq n})_{\alpha \in\Natural^k}
\stackrel{d}{=} ( w_\alpha, (g^\prime_{\xi',i}(\alpha))_{i\leq n}
)_{\alpha\in\Natural^k}
\end{equation}
have the same distribution, where $(g^\prime_{\xi',i}(\alpha))$ is a
random field (no longer Gaussian)
associated with the Ruelle probability cascades defined from the
Gaussian field $(g_{\xi',i}(\alpha))$
by an explicit change of density; see equation (7) in~\cite{PTDR}.
Therefore, (\ref{Invarg})
can be rewritten as
%
%e1.70 ###
\begin{equation}\label{InvarDR2}
\mathbb{E}\prod_{l\leq q}\mathbb{E}' \prod_{i\in C_l} \bar{\sigma}_i
=\mathbb{E}\prod_{l\leq q} \sum_{\alpha} w_\alpha\prod_{i\in C_l}
\th  g^\prime_{\xi',i}(\alpha),
\end{equation}
which can now be interpreted as the explicit construction of $\bar
{\sigma}(w,u,v)$.
The first coordinate $w$ corresponds to generating the weights
$(w_\alpha)_{\alpha\in\Natural^k}$
of the Ruelle probability cascade with the parameters
$0=m_1<\cdots<m_k<1$,
the second coordinate $u$ plays the role of sampling an index $\alpha$
according to the weights
$(w_\alpha)$ and the last coordinate $v$ corresponds to generating the
random field
$(g_{\xi'}^\prime(\alpha))$, so that the directing measure $\eta_w$
carries weight $w_\alpha$
at the point $\th  g^\prime_{\xi'}(\alpha)$ in $\LL^2([0,1],dv)$.
Another way to write this is to consider
a partition $(C_\alpha)_{\alpha\in\Natural^k}$ of $[0,1]$ into
intervals of length
$|C_\alpha|=w_\alpha$ and let
%
%e1.71 ###
\begin{equation}\label{orderSK}
\bar{\sigma}(w,u,v) =
\sum_{\alpha\in\Natural^k} I(u\in C_\alpha)  \th  g^\prime
_{\xi'}(\alpha),
\end{equation}
where we\vspace*{2pt} keep the dependence of $(C_\alpha)$ on $w$ and $(g^\prime
_{\xi'}(\alpha))$ on $v$
implicit. In particular, (\ref{InvarDR2}) implies that the limiting
distribution of the Gibbs\vadjust{\goodbreak}
averages $\la\sigma_i\ra$ of finitely many spins $1\leq i\leq n$
coincides with the distribution of
%
%e1.72 ###
\begin{equation}
\sum_{\alpha} w_\alpha\th  g^\prime_{\xi',i}(\alpha)
\qquad\mbox{for $1\leq i\leq n$.}
\end{equation}
This can be thought of as the generalization of the high temperature
result (Theorem 2.4.12
in~\cite{SG}) under the assumption of the Parisi ultrametricity. It
will be clear from the proof of
Theorem~\ref{ThSCSK} that the right-hand side of (\ref{Invarg}) is
continuous with respect to
the distribution of the overlap array (\ref{overlaps}) and, on the
other hand, it is well known
that ultrametricity allows one to approximate any overlap array by a
discretized overlap
array satisfying (\ref{finite}) uniformly while preserving
ultrametricity and the Ghirlanda--Guerra
identities. Therefore, one can think of the case of an arbitrary
distribution of the overlap simply
as the limiting case of the above construction for discrete overlaps.

\subsubsection*{Parisi ansatz in the diluted models}
To make a transition to the case of diluted models let us look more
closely at equation
(\ref{orderSK}). Original Gaussian field $(g_{\xi'}(\alpha))$
indexed by
$\alpha=(\alpha_1,\ldots, \alpha_k)\in\Natural^k$ associated to
the Ruelle probability cascades
is of the form~\cite{AA}
\[
g_{\xi'}(\alpha) = g_{\xi'}(\alpha_1)+g_{\xi'}(\alpha_1,\alpha
_2)+\cdots+
g_{\xi'}(\alpha_1,\ldots,\alpha_k),
\]
where random variables $g_{\xi'}(\alpha_1,\ldots,\alpha_l)$ are
Gaussian with variances
$\xi'(q_l)-\xi'(q_{l-1})$ independent for different $1\leq l\leq k$
and different
$(\alpha_1,\ldots,\alpha_l)$. The field $(g_{\xi'}^\prime(\alpha
))$ on the right-hand side of
(\ref{gprime}) is again of the form
\[
g_{\xi'}^\prime(\alpha) = g_{\xi'}^\prime(\alpha_1)+g_{\xi
'}^\prime(\alpha_1,\alpha_2)+\cdots+
g_{\xi'}^\prime(\alpha_1,\ldots,\alpha_k),
\]
and for each $l\leq k$ the sequence $(g_{\xi'}^\prime(\alpha
_1,\ldots,\alpha_l))_{\alpha_l\geq1}$
is i.i.d. from distribution defined by the explicit change of density
(equation (7) in~\cite{PTDR})
which depends on $g_{\xi'}^\prime(\alpha_1), \ldots, g_{\xi
'}^\prime(\alpha_1,\ldots,\alpha_{l-1})$,
and these sequences are independent for different $(\alpha_1,\ldots,
\alpha_{l-1})$ conditionally
on the sequences $(g_{\xi'}^\prime(\alpha_1)), \ldots, (g_{\xi
'}^\prime(\alpha_1,\ldots,\break\alpha_{l-1}))$.
This means that one can generate the process $(g_{\xi'}^\prime(\alpha
))$ recursively as follows.
Let $v(\alpha_1,\ldots,\alpha_l)$ be random variables uniform on
$[0,1]$ independent for different
$1\leq l\leq k$ and different $(\alpha_1,\ldots,\alpha_l)$. Then for
$1\leq l\leq k$ we can define
%
%e1.73 ###
\begin{equation}\label{Ql}\qquad
g_{\xi'}^\prime(\alpha_1,\ldots,\alpha_l)
=Q_l(g_{\xi'}^\prime(\alpha_1), \ldots, g_{\xi'}^\prime(\alpha
_1,\ldots,\alpha_{l-1}),
v(\alpha_1,\ldots,\alpha_l)),
\end{equation}
where $Q_l$ as a function of the last variable is the quantile
transform of the distribution defined
by the aforementioned change of density. Combining all the steps of the
recursion we get
%
%e1.74 ###
\begin{equation}\label{QSK}
g_{\xi'}^\prime(\alpha) = Q(v(\alpha_1),\ldots,v(\alpha_1,\ldots
,\alpha_k))
\end{equation}
for some specific function $Q$. Equation (\ref{orderSK}) becomes
%
%e1.75 ###
\begin{equation}\label{orderSK2}
\bar{\sigma}(w,u,v) = \sum_{\alpha\in\Natural^k} I(u\in C_\alpha)
\varphi(v(\alpha_1),\ldots,v(\alpha_1,\ldots,\alpha_k)),\vadjust{\goodbreak}
\end{equation}
where $\varphi= \th\circ Q$, and again, as in (\ref
{orderSK}), we keep the dependence
of $(C_\alpha)$ on $w$ and $(v(\alpha_1,\ldots,\alpha_l))$ on $v$ implicit.
Let us emphasize that the change of density that defines $Q_l$ in (\ref
{Ql}) and, therefore,
the functions $Q, \varphi$ and $\bar{\sigma}$ are completely
determined by the parameters of
the distribution of one overlap in (\ref{finite}) which is the
functional order parameter of the Parisi
ansatz in the Sherrington--Kirkpatrick model.
What seems to be the main (and only) difference in the Parisi ansatz
for diluted models is that
this function $\varphi$ is allowed to be an arbitrary $(-1,1)$ valued
function, which we will now explain.

The Parisi functional order parameter in the diluted models appears in
the description of the free
energy, and one can make the connection to the generic functional order
parameter $\bar{\sigma}$
by comparing the Parisi formula for the free energy to the
representation (\ref{CalP}), (\ref{FE}).
For example, in the notation of~\cite{PT} where the order parameter
was encoded by the Ruelle
probability cascade weights $(w_\alpha)$ and associated random field
$(x({\alpha}))$
for $\alpha\in\Natural^k$, it is easy to see that in order for (\ref
{CalP}) to match the Parisi formula in
\cite{PT}, $\bar{\sigma}$ should be defined exactly as
in (\ref{orderSK}),
%
%e1.76 ###
\begin{equation}\label{orderdiluted}
\bar{\sigma}(w,u,v) =
\sum_{\alpha\in\Natural^k} I(u\in C_\alpha)  \th  x(\alpha).
\end{equation}
The only difference from (\ref{orderSK}) is how the random field
$(x(\alpha))$ is generated
compared to $(g_{\xi'}^\prime(\alpha))$, and once we recall how
$(x(\alpha))$ is generated
according to the Parisi ansatz, we will realize that one can write
exactly the same representation as~(\ref{QSK}),
%
%e1.77 ###
\begin{equation}\label{Qdiluted}
x(\alpha) = Q(v(\alpha_1),\ldots,v(\alpha_1,\ldots,\alpha_k)),
\end{equation}
only now $Q$ is allowed to be arbitrary. The field $(x(\alpha))$ is
customarily generated as follows.
Let ${\mathrm P}_1$ be the set of probability measures on $\Reals$, and by
induction on $l\leq k$ we define
${\mathrm P}_{l+1}$ as the set of probability measures on~${\mathrm P}_{l}$.
Let us fix $\eta\in{\mathrm P}_{k}$
(the basic parameter) and define a random sequence
$(\eta(\alpha_1), \ldots, \eta(\alpha_1, \ldots, \alpha_{k-1}),
x(\alpha_1,\ldots, \alpha_k))$ as follows.
Given $\eta$, the sequence $(\eta(\alpha_1))_{\alpha_1 \geq1}$ of
elements of ${\mathrm P}_{k-1}$ is
i.i.d. from distribution $\eta$. For $1\leq l\leq k-1$, given all the elements
$\eta(\alpha_1, \ldots, \alpha_{s})$ for all values of the integers
$\alpha_1, \ldots, \alpha_s$
and all $s \leq l-1$, the sequence $(\eta(\alpha_1, \ldots, \alpha
_l))_{\alpha_l \geq1}$ of
elements of ${\mathrm P}_{k-l}$ is i.i.d. from distribution $\eta(\alpha
_1, \ldots, \alpha_{l-1})$,
and these sequences are independent of each other for different values of
$(\alpha_1, \ldots, \alpha_{l-1})$. Finally, given all the elements
$\eta(\alpha_1, \ldots, \alpha_{s})$
for all values of the integers $\alpha_1, \ldots, \alpha_s$ and all
$s \leq k-1$ the sequence
$(x(\alpha_1,\ldots,\alpha_k))_{\alpha_k\geq1}$ is i.i.d. on
$\Reals$ with distribution
$\eta(\alpha_1,\ldots,\alpha_{k-1})$ and these sequences are
independent for different values
of $(\alpha_1,\ldots,\alpha_{k-1})$. The process of generating $x$'s
can be represented
schematically as
%
%e1.78 ###
\begin{equation}\label{xs}
\eta\to\eta(\alpha_1)\to\cdots\to\eta(\alpha_1,\ldots,\alpha
_{k-1})\to
x(\alpha_1,\ldots,\alpha_k).
\end{equation}
Now, as above, let $v(\alpha_1,\ldots,\alpha_l)$ be random variables
uniform on $[0,1]$
independent for different $1\leq l\leq k$ and different $(\alpha
_1,\ldots,\alpha_l)$.
First, random variables $(\eta(\alpha_1))_{\alpha_1 \geq1}$ are
i.i.d. from probability
measure $\eta$ on ${\mathrm P}_{k-1}$\vadjust{\goodbreak} and, therefore, can be generated as
%
%e1.79 ###
\begin{equation}\label{oned}
\eta(\alpha_1) = Q_{k-1}(v(\alpha_1))
\end{equation}
for some function $Q_{k-1}\dvtx[0,1]\to{\mathrm P}_{k-1}$. Next, random variables
$(\eta(\alpha_1,\alpha_2))_{\alpha_2 \geq1}$ are i.i.d. from probability
measure $\eta(\alpha_1)$ on ${\mathrm P}_{k-2}$ and, therefore, can be
generated as
\[
\eta(\alpha_1,\alpha_2) = \tilde{Q}_{k-2}(\eta(\alpha_1),v(\alpha
_1,\alpha_2))
\]
for some function $ \tilde{Q}_{k-2}(\eta(\alpha_1),\cdot)\dvtx[0,1]\to
{\mathrm P}_{k-2}$.
Combining with (\ref{oned}), we can write
%
%e1.80 ###
\begin{equation}\label{twod}
\eta(\alpha_1,\alpha_2) = Q_{k-2}(v(\alpha_1),v(\alpha_1,\alpha_2))
\end{equation}
for some function $Q_{k-2}\dvtx[0,1]^2\to{\mathrm P}_{k-2}$. We can continue
this construction
recursively and at the end we will get
%
%e1.81 ###
\begin{equation}\label{kd}
x(\alpha_1,\ldots,\alpha_k) = Q(v(\alpha_1),\ldots,v(\alpha
_1,\ldots,\alpha_k))
\end{equation}
for some function $Q\dvtx[0,1]^k\to\Reals$, which is exactly (\ref{Qdiluted}).
This representation gives some choice of $Q$ for a given $\eta\in{\mathrm P}_k$, but any choice of $Q$
corresponds to some $\eta$, which is obvious by reverse induction and
identifying a function
of uniform r.v. on $[0,1]$ with the distribution on its image.

To summarize, the Parisi ansatz can be
expressed in terms of $\bar{\sigma}$ by saying that equation (\ref
{orderSK2}) must hold
for some choice of $(-1,1)$ valued function $\varphi$. Of course, in
general this statement should
be understood in the limiting sense when the number $(k-1)$ of
replica-symmetry breaking steps goes
to infinity. Precise statement should be that in the diluted models any
limiting distribution $\mu\in\M$
of the array (\ref{sigma}) over a subsequence of $(\mu_N)$ can be
approximated by the distribution
of the array generated by $\bar{\sigma}(w,u,v)$ as in (\ref
{orderSK2}) for large enough $k$,
some function $\varphi\dvtx[0,1]^k\to(-1,1)$ and some parameters
$0=m_1<\cdots<m_k<1$
of the distribution of weights $(w_\alpha)$ in the Ruelle probability cascades.

This formulation clarifies another statement of the physicists, namely,
that multi-overlap
$R^\infty_{1,\ldots, n}$ in (\ref{multioverlap}) is the function of
the overlaps
$R^\infty_{l,l'}$ in (\ref{overlaps}) for $1\leq l<l'\leq n$.
According to (\ref{orderSK2})
the choice of $u_1,\ldots,u_n$ corresponds to the choice of indices
$\alpha^1,\ldots,\alpha^n\in\Natural^k$ so that
\[
R^\infty_{1,\ldots, n}
=
\mathbb{E}\varphi(v(\alpha_1^1),\ldots,v(\alpha_1^1,\ldots,\alpha_k^1))
\cdots
\varphi(v(\alpha_1^n),\ldots,v(\alpha_1^n,\ldots,\alpha_k^n)).
\]
On the other hand, if we denote $\alpha^1\wedge\alpha^2 = \min\{i\dvtx
\alpha_i^1\not= \alpha_i^2\}$
and $\alpha^1\wedge\alpha^2=k+1$ if $\alpha^1=\alpha^2$, then the
overlap takes finitely
many values
\begin{eqnarray*}
R^\infty_{1,2}
&=&
\mathbb{E}\varphi(v(\alpha_1^1),\ldots,v(\alpha_1^1,\ldots,\alpha_k^1))
\varphi(v(\alpha_1^2),\ldots,v(\alpha_1^2,\ldots,\alpha_k^2))\\
&=&
q_{\alpha^1\wedge\alpha^2}
\end{eqnarray*}
for some $0\leq q_1\leq\cdots\leq q_{k+1}\leq1$. This means\vspace*{1pt} that the
values of the overlaps
$(R^\infty_{l,l'})$ determine $(\alpha^l\wedge\alpha^{l'})$ for
$1\leq l<l'\leq n$. It is also clear
that the multi-overlap $R^\infty_{1,\ldots, n}$ is the same for two
sets of indices
$(\alpha^1,\ldots,\alpha^n)$ and $(\beta^1,\ldots,\beta^n)$ for which
$(\alpha^l\wedge\alpha^{l'})=(\beta^{\rho(l)}\wedge\beta^{\rho
(l')})$ for some permutation
$\rho$ the set $\{1,\ldots, n\}$. In this sense, given representation
(\ref{orderSK2}), the overlaps
indeed determine the value of the multi-overlap.
At the moment we have no idea how (\ref{orderSK2}) can be proved, but
it is helpful to have
a point of view that formulates precisely the predictions of the Parisi ansatz.

While many technical details will be quite different, the main line of
the arguments in the setting
of the Sherrington--Kirkpatrick model in Section~\ref{SecSKK} will be
parallel to the arguments
in Section~\ref{SecDiluted} for diluted models. A reader only
interested in the Sherrington--Kirkpatrick model should read Lemma \ref
{LemShift} before skipping to Section~\ref{SecSKK}.

%s2 ###
\section{Diluted models}\label{SecDiluted}

%s2.1 ###
\subsection{Properties of convergence}\label{SecPC}

Let us first record a simple consequence of the fact that the
distribution of the array in (\ref{sigma})
is the limit of the distribution of spins $(\sigma_i^l)$ under the
annealed product Gibbs measure.
As usual, $\la\cdot\ra$ will denote the expectation with respect to
the random Gibbs measure.
Also, recall the definition of $\mathbb{E}'$ before Theorem~\ref{ThSC}.
\begin{lemma}\label{LemP1}
Let $h_1,\ldots, h_m\dvtx\{-1,+1\}^n\to[-K,K]$ be some bounded functions
of $n$ spins, and let
$h$ be a continuous function on $[-K,K]^m$. Let $\vsi= (\sigma
_i)_{1\leq i\leq n}$,
and let $\mathbf{s}=(s_i^1)_{1\leq i\leq n}$ defined in (\ref{sigma})
for some $\mu\in\M$.
If $\mu_N$ converges to $\mu$ over subsequence $(N_k)$, then
%
%e2.1 ###
\begin{equation}\label{moments}
\lim_{N_k\to\infty} \mathbb{E}h(\la h_1(\vsi)\ra,\ldots, \la
h_m(\vsi)\ra) = \mathbb{E}h( \mathbb{E}' h_1(\mathbf{s}),\ldots,
\mathbb{E}' h_m(\mathbf{s})).
\end{equation}
\end{lemma}
\begin{pf}
Since it is enough to prove this for polynomials $h$ and since each
$h_l$ is a polynomial in its
coordinates, this statement is simply a convergence of moments
\[
\lim_{N_k\to\infty} \mathbb{E}\Bigl\la\prod\sigma_i^l\Bigr\ra
= \mathbb{E}\prod s_i^l,
\]
where the product is over a finite subset of indices $(i,l)$.
\end{pf}

We will often use this lemma for random functions $h, (h_l)$
independent of all other randomness,
simply by applying (\ref{moments}) conditionally on the randomness of
these functions. Justifications
of convergence will always be omitted because of their triviality.

Another simple property of convergence of spin distributions under the
annealed Gibbs measure
in diluted models is that adding or removing a finite number of terms
to the Poisson number of terms
$\pi(\alpha N)$ or $\pi(c_N)$ in (\ref{Hampluspert}) does not affect\vadjust{\goodbreak}
the limit of these distribution over any subsequence for which the
limit exists. Let $(N_k)_{k\geq1}$ be any such subsequence, and let
$n, m$ be fixed integers.
In fact, it will be clear from the proof that one can let $n,m$ grow
with $N_k$,
but we will not need this. Let $H_{N}'$ be defined exactly as (\ref
{Hampluspert})
only with $\pi(\alpha N) + n$ terms instead of $\pi(\alpha N)$ in the
first sum and
$\pi(c_N)+m$ instead of $\pi(c_N)$ in the perturbation term, and
let $\la\cdot\ra'$ denote the corresponding Gibbs measure.
\begin{lemma}\label{LemShift}
For any bounded function $h$ of finitely many spins in array $(\sigma
_i^l)$ we have
%
%e2.2 ###
\begin{equation}\label{Shift}
{\lim_{N\to\infty}} | \mathbb{E}\la h \ra' -\mathbb{E}\la h\ra | =0.
\end{equation}
\end{lemma}
\begin{pf}
For certainty, let us assume that $n,m\geq0$ and
$|h|\leq1$. If we denote by
$\la\cdot\ra_{i,j}$ the Gibbs average conditionally on $\pi(\alpha
N) = i$ and $\pi(c_N)=j$, then
\[
\mathbb{E}\la h\ra =\sum_{i,j\geq0} \pi(\alpha N,i) \pi
(c_N,j)\mathbb{E}\la h\ra_{i,j},
\]
where from now on $\pi(\lambda,k)=\lambda^k e^{-\lambda}/k!$ and
\begin{eqnarray*}
\mathbb{E}\la h\ra'
&=&
\sum_{i,j\geq0} \pi(\alpha N,i)\pi(c_N,j) \mathbb{E}\la h\ra_{i+n,j+m}
\\
&=&
\sum_{i\geq n,j\geq m} \pi(\alpha N,i-n)\pi(c_N,j-m) \mathbb{E}\la
h\ra_{i,j}.
\end{eqnarray*}
Therefore,
\begin{eqnarray*}
| \mathbb{E}\la h \ra' -\mathbb{E}\la h\ra |
&\leq& \sum_{i<n}
\pi(\alpha N,i) + \sum_{j<m} \pi(c_N,j)
\\
&&{}+
\sum_{i\geq n,j\geq m} | \pi(\alpha N,i-n)\pi(c_N,j-m) -\pi(\alpha
N,i)\pi(c_N,j)|
\\
&\leq&
\sum_{i<n} \pi(\alpha N,i)
+
\sum_{j<m} \pi(c_N,j)
+
\sum_{i\geq n} | \pi(\alpha N,i-n) -\pi(\alpha N,i)|\\
&&{}+
\sum_{j\geq m} | \pi(c_N,j-m) -\pi(c_N,j)|.
\end{eqnarray*}
The first two sums obviously go to zero. One can see that the third sum
goes to zero as follows.
Poisson distribution with mean $\alpha N$ is concentrated inside the range
%
%e2.3 ###
\begin{equation}\label{irange}
\alpha N - \sqrt{N\log N} \leq i \leq\alpha N+\sqrt{N\log N}.
\end{equation}
If we write
%
%e2.4 ###
\begin{equation}\label{nappr}
| \pi(\alpha N,i-n) -\pi(\alpha N,i)| = \pi(\alpha N,i)\biggl| 1-
\frac{i!}{(i-n)!} (\alpha N)^{-n} \biggr|,\vadjust{\goodbreak}
\end{equation}
then it remains to note that\vspace*{-1pt}
\[
\frac{i!}{(i-n)!} (\alpha N)^{-n}
=
\frac{i(i-1)\cdots(i-n+1)}{(\alpha N)^n}
\to1\vspace*{-1pt}
\]
uniformly inside the range (\ref{irange}). Similarly, the last sum
goes to zero which finishes the proof.\vspace*{-2pt}
\end{pf}
\begin{Remark*}
Lemma~\ref{LemShift} implies that (\ref{Shift})
holds even if $n$ is a random variable.
We will use this observation in the case when $H_{N}'$ is defined
exactly as
(\ref{Hampluspert}) only with $\pi(\alpha N + n)$ terms instead of
$\pi(\alpha N)$.
In fact, in this case one can write\vspace*{-1pt}
\[
\mathbb{E}\la h\ra'
=
\sum_{i,j\geq0} \pi(\alpha N+n,i)\pi(c_N,j) \mathbb{E}\la h\ra_{i,j+m}\vspace*{-1pt}
\]
and instead of (\ref{nappr}) use\vspace*{-1pt}
\[
| \pi(\alpha N+n,i) -\pi(\alpha N,i)| = \pi(\alpha N,i)\biggl| 1-
\biggl(1+\frac{n}{\alpha N}\biggr)^i e^{-n} \biggr|\vspace*{-1pt}
\]
and notice that again the last factor goes to zero uniformly over range
(\ref{irange}).
Similarly, one can have $\pi(c_N + n)$ instead of $\pi(c_N)$ terms in
the perturbation
Hamiltonian without affecting convergence.\vspace*{-2pt}
\end{Remark*}

Due to the perturbation term (\ref{Hampert}) the following important
property of convergence holds.
\begin{lemma}\label{LemShift2}
If $\mu_N$ converges to $\mu$ over subsequence $(N_k)$ then it also converges
to $\mu$ over subsequence $(N_k+n)$ for any $n\geq1$.\vspace*{-2pt}
\end{lemma}
\begin{pf}
We will show that the joint moments of spins converge
to the same limit over
subsequences that differ by a finite shift $n$. Let $h=\prod_{j\leq q}
h_j$ where
$h_j = \prod_{i\in C_j} \sigma_i^j$ over some finite sets of spin
coordinates $C_j$.
Let us denote by $\la\cdot\ra_N$ the
Gibbs average with respect to the Hamiltonian (\ref{Hampluspert})
defined on $N$ coordinates.
We will show that\vspace*{-1pt}
\[
{\lim_{N\to\infty}} |\mathbb{E}\la h\ra_{N+n} - \mathbb{E}\la h\ra
_N| = 0.\vspace*{-1pt}
\]
Let us rewrite $\mathbb{E}\la h\ra_{N+n}$ by treating the last $n$
coordinates as cavity coordinates.
Let us separate the $\pi(\alpha(N+n))$ terms in the first sum\vspace*{-1pt}
%
%e2.5 ###
\begin{equation}\label{firstsum}
\sum_{k\leq\pi(\alpha(N+n))}\theta_k(\sigma_{i_{1,k}}, \ldots
,\sigma_{i_{p,k}})\vspace*{-1pt}
\end{equation}
of the Hamiltonian $H_{N+n}(\vsi)$ in (\ref{Hampluspert}) into
several groups:

\begin{longlist}[($2l$)]
\item[(1)]
terms for $k$ such that all indices $i_{1,k},\ldots, i_{p,k}\leq
N$;

For $1\leq l\leq n$:

\item[($2l$)] terms with exactly one of indices $i_{1,k},\ldots, i_{p,k}$
equal to $N+l$ and all others~$\leq$~$N$;

\item[(3)] terms with at least two of indices $i_{1,k},\ldots, i_{p,k}\geq
N$.\vadjust{\goodbreak}
\end{longlist}
The probabilities that a term is of these three type are
\[
p_1=\biggl(\frac{N}{N+n} \biggr)^p,\qquad
p_{2,l} = p\frac{1}{N+n} \biggl(\frac{N}{N+n}\biggr)^{p-1},\qquad
p_3=1-p_1-\sum_{l\leq n}p_{2,l}.
\]
Therefore, the number of terms in these groups are independent
Poisson random variables with means
\begin{eqnarray*}
\alpha(N+n) p_1 &=& \alpha(N +n- np) + O(N^{-1}),\\
\alpha(N+n) p_{2,l} &=& \alpha p + O(N^{-1}),\\
\alpha(N+n) p_3 &=& O(N^{-1}).
\end{eqnarray*}
We can redefine the number of terms in each group to be exactly of
means $\alpha(N+n-np),
\alpha p$ and $0$ since asymptotically it does not affect $\mathbb
{E}\la h\ra_{N+n}$ as in Lem\-ma~\ref{LemShift}
or using assumption (\ref{condition4}).
Thus, if we write $\vsi=(\bolds{\rho},\bolds{\eps})\in\Sigma_{N+n}$
for the first
$N$ coordinates $\bolds{\rho}=(\rho_1,\ldots,\rho_{N})$ and the last $n$
cavity coordinates $\bolds{\eps}=(\eps_1,\ldots, \eps_{n})$, then
(\ref{firstsum})
can be replaced with
%
%e2.6 ###
\begin{eqnarray}\label{firstterm}
&&\sum_{k\leq\pi(\alpha(N+n-np))}\theta_k(\rho_{i_{1,k}},
\ldots,\rho_{i_{p,k}})\nonumber\\[-8pt]\\[-8pt]
&&\qquad{}+\sum_{l\leq n}
\sum_{k\leq\pi_l(\alpha p)}\theta_{k,l}(\eps_l, \rho_{i_{1,k,l}},
\ldots,\rho_{i_{p-1,k,l}}),\nonumber
\end{eqnarray}
where indices $i_{1,k},\ldots, i_{p,k}$ and $i_{1,k,l},\ldots,
i_{p-1,k,l}$ are all
uniformly distributed on $\{1,\ldots, N\}$. Let us now consider the
perturbation term
in (\ref{Hampluspert})
%
%e2.7 ###
\begin{equation}\label{pertn}
\sum_{l\leq\pi(c_{N+n})} \log\Av_\eps \exp\sum_{k\leq\hat{\pi}_l(\alpha
p)} \hat{\theta}_{k,l}(\eps, \sigma_{j_{1,k,l}},
\ldots,\sigma_{j_{p-1,k,l}}),
\end{equation}
where $j_{1,k,l},\ldots, j_{p-1,k,l}$ are uniformly distributed on $\{
1,\ldots, N+n\}$.
Here, we used independent copies $\hat{\pi}_l$ and $\hat{\theta
}_{k,l}$ since
$\pi_l$ and $\theta_{k,l}$ were already used in (\ref{firstterm}).
The expected number of all such indices in (\ref{pertn}) that belong
to $\{N+1,\ldots, N+n\}$ is
$c_{N+n}\alpha p(p-1)n/(N+n)\to0$ which means that with high probability
all indices belong to $\{1,\ldots, N\}$. As a result, asymptotically
$\mathbb{E}\la h\ra_{N+n}$
will not be affected if we replace the perturbation term (\ref{pertn}) with
%
%e2.8 ###
\begin{equation}\label{pertn2}
\sum_{l\leq\pi(c_{N+n})} \log\Av_\eps \exp\sum_{k\leq\hat{\pi}_l(\alpha
p)} \hat{\theta}_{k,l}(\eps, \rho_{j_{1,k,l}},
\ldots,\rho_{j_{p-1,k,l}}),
\end{equation}
where $j_{1,k,l},\ldots, j_{p-1,k,l}$ are uniformly distributed on $\{
1,\ldots, N\}$.
Thus, we can assume from now on that $\mathbb{E}\la h\ra_{N+n}$ is
computed with respect to the Hamiltonian
which is the sum of (\ref{firstterm}) and (\ref{pertn2}). If $\la
\cdot\ra_N'$ denotes the Gibbs average
on $\Sigma_N$ with respect to the Hamiltonian
\begin{eqnarray*}
-H_N'(\bolds{\rho})
&=&
\sum_{k\leq\pi(\alpha(N+n-np))}\theta_k(\rho_{i_{1,k}},
\ldots,\rho_{i_{p,k}})
\\
&&{}
+
\sum_{l\leq\pi(c_{N+n})}
\log\Av_\eps
\exp\sum_{k\leq\hat{\pi}_l(\alpha p)}
\hat{\theta}_{k,l}(\eps, \rho_{j_{1,k,l}}, \ldots,\rho_{j_{p-1,k,l}}),
\end{eqnarray*}
then each factor in
\[
\la h\ra_{N+n}=\prod_{j\leq q} \la h_j \ra_{N+n}
=\prod_{j\leq q} \biggl\la\prod_{i\in C_j}\sigma_i \biggr\ra_{N+n}
=\prod_{j\leq q} \biggl\la\prod_{i\in C_j}\rho_i \biggr\ra_{N+n}
\]
can be written as
\begin{eqnarray*}
\la h_j \ra_{N+n}
&=&
\frac{
\la\prod_{i\in C_j}\rho_i
\Av_\eps\exp\sum_{l\leq n}
\sum_{k\leq\pi_l(\alpha p)}\theta_{k,l}(\eps_l, \rho_{i_{1,k,l}},
\ldots,\rho_{i_{p-1,k,l}})\ra_N'}{
\la
\Av_\eps\exp\sum_{l\leq n}
\sum_{k\leq\pi_l(\alpha p)}\theta_{k,l}(\eps_l, \rho_{i_{1,k,l}},
\ldots,\rho_{i_{p-1,k,l}})\ra_N'}
\\
&=&
\biggl\la\prod_{i\in C_j}\rho_i \biggr\ra_N'',
\end{eqnarray*}
where $\la\cdot\ra_N''$ is the Gibbs average on $\Sigma_N$
corresponding to
the Hamiltonian
\[
-H_N''(\bolds{\rho})
=
-H_N'(\bolds{\rho})
+
\sum_{l\leq n} \log\Av_\eps
\exp
\sum_{k\leq\pi_l(\alpha p)}\theta_{k,l}(\eps, \rho_{i_{1,k,l}},
\ldots,\rho_{i_{p-1,k,l}}).
\]
But this Hamiltonian differs from the original Hamiltonian (\ref
{Hampluspert}) only
in that the first sum has $\pi(\alpha(N+n-np))$ terms instead of $\pi
(\alpha N)$,
and the perturbation term has $\pi(c_{N+n})+n$ terms instead of
$\pi(c_N)$.
Therefore, appealing to Lemma~\ref{LemShift} and remark after it shows that
$\mathbb{E}\la h\ra_N''$ is asymptotically equivalent to $\mathbb
{E}\la h\ra_N$ and this finishes the proof.
\end{pf}

%s2.2 ###
\subsection{Lower bound}\label{SecLB}

\begin{lemma}\label{LemLB}
There exists $\mu\in\M$ such that
$
\lim_{N\to\infty} F_N \geq\P(\mu).
$
\end{lemma}
\begin{pf}
We will obtain the lower bound using the well-known fact that
%
%e2.9 ###
\begin{equation}\label{Zdiff}
\lim_{N\to\infty} F_N \geq \liminf_{N\to\infty} \bigl((N+1)F_{N+1} - NF_N\bigr)
=\liminf_{N\to\infty} \mathbb{E}\log\frac{Z_{N+1}}{Z_N}.
\end{equation}
Suppose that this lower limit is achieved over subsequence $(N_k)$,
and let $\mu\in\M$ be a limit of $(\mu_N)$ over some subsubsequence
of $(N_k)$.
Let $\sigma=\sigma_\mu$. The considerations will be very similar to
the proof
of Lemma~\ref{LemShift2}. Let us consider $\mathbb{E}\log Z_{N+1}$,
and let us start by
separating the $\pi(\alpha(N+1))$ terms in the first sum in the
Hamiltonian $H_{N+1}$
in (\ref{Hampluspert}) into three groups:
(1) terms for $k$ such that all indices $i_{1,k},\ldots, i_{p,k}\leq N$;
(2) terms with exactly one of indices $i_{1,k},\ldots, i_{p,k}$ equal
to $N+1$;
(3) terms with at least two of indices $i_{1,k},\ldots, i_{p,k}$
equal\vadjust{\goodbreak}
to $N+1$.
The probabilities that a term is of these three types are
\[
p_1=\biggl(\frac{N}{N+1}\biggr)^p,\qquad
p_2 = p\frac{1}{N+1} \biggl(\frac{N}{N+1}\biggr)^{p-1},\qquad
p_3=1-p_1-p_2
\]
correspondingly. Therefore, the number of terms in these three groups
are independent
Poisson random variables with means
\begin{eqnarray*}
\alpha(N+1) p_1 &=& \alpha(N -p+1) + O(N^{-1}),\\
\alpha(N+1) p_2 &=& \alpha p + O(N^{-1}),\\
\alpha(N+1) p_3 &=& O(N^{-1}).
\end{eqnarray*}
For simplicity of notation, let us pretend that the number of terms in
each group is exactly
of means $\alpha(N-p), \alpha p$ and $0$ since it will be clear from
considerations below that asymptotically it does not affect the limit
in (\ref{Zdiff}).
If we write $\vsi=(\bolds{\rho},\eps)\in\Sigma_{N+1}$ for $
\bolds{\rho}\in\Sigma_N$ and $\eps\in\{-1,+1\}$,
then we can write the first term in $H_{N+1}(\vsi)$ as
%
%e2.10 ###
\begin{equation}\label{Hamprime}\qquad
\sum_{k\leq\pi(\alpha(N-p+1))}\theta_k(\rho_{i_{1,k}},\ldots
,\rho_{i_{p,k}}) +\sum_{k\leq\pi(\alpha p)}\hat{\theta}_k(\eps,
\rho_{j_{1,k}}, \ldots,\rho_{j_{p-1,k}}),
\end{equation}
where indices $i_{1,k},\ldots, i_{p,k}$ and $j_{1,k},\ldots,
j_{p-1,k}$ are uniformly distributed
on $\{1,\ldots, N\}$.
Similarly, we could split the $\pi(c_{N+1})$ terms in the perturbation
Hamiltonian
(\ref{Hampert}) into indices $l$ for which all $i_{1,k,l},\ldots
,i_{p-1,k,l}\leq N$
and indices $l$ for which at least one of these indices equals $N+1$.
However, as in the proof of Lemma~\ref{LemShift2}, since with high
probability all these
indices will be $\leq N$ and $|c_{N+1}-c_N|\to0$,
we can simply replace the perturbation term with
%
%e2.11 ###
\begin{equation}\label{pertn3}
\sum_{l\leq\pi(c_{N})} \log\Av_\eps \exp\sum_{k\leq\pi_l(\alpha p)}
\theta_{k,l}(\eps, \rho_{i_{1,k,l}}, \ldots,\rho_{i_{p-1,k,l}}),
\end{equation}
where $i_{1,k,l},\ldots, i_{p-1,k,l}$ are uniformly distributed on $\{
1,\ldots, N\}$.
Let $\la\cdot\ra'$ be the Gibbs average on $\Sigma_N$ corresponding
to the
Hamiltonian
\begin{eqnarray*}
-H_N'(\bolds{\rho})
&=&
\sum_{k\leq\pi(\alpha(N-p+1))}\theta_k(\rho_{i_{1,k}},\ldots
,\rho_{i_{p,k}})
\\
&&{}
+\sum_{l\leq\pi(c_{N})}
\log\Av_\eps
\exp\sum_{k\leq\pi_l(\alpha p)}
\theta_{k,l}(\eps, \rho_{i_{1,k,l}}, \ldots,\rho_{i_{p-1,k,l}})
\end{eqnarray*}
and $Z_N'$ be the corresponding partition function.
Then
%
%e2.12 ###
\begin{equation}\label{this}\quad
\mathbb{E}\log\frac{Z_{N+1}}{Z_N'}
=\mathbb{E}\log\biggl\la\sum_{\eps=\pm1} \exp
\sum_{k\leq\pi(\alpha p)}\hat{\theta}_k(\eps, \rho_{j_{1,k}},
\ldots,\rho_{j_{p-1,k}})
\biggr\ra'.
\end{equation}
Conditionally on $\pi(\alpha p)$ and $(\hat{\theta}_k)$ and on the
event that all
indices $j_{1,k},\ldots,\break j_{p-1,k}$ are different, Lemmas~\ref{LemP1}
and~\ref{LemShift}
imply that (\ref{this}) converges to
\[
\mathbb{E}\log\mathbb{E}' \sum_{\eps=\pm1} \exp\sum_{k\leq\pi
(\alpha p)}
\hat{\theta}_k(\eps, s_{1,k},\ldots, s_{p-1,k}).
\]
For large $N$, indices $j_{1,k},\ldots, j_{p-1,k}$ will all be different
for all $k\leq\pi(\alpha p)$ with high probability and, therefore,
this convergence
holds unconditionally. Similarly, one can analyze $\mathbb{E}\log
(Z_{N}/Z_N')$.
Let us split the first sum in the definition of $-H_N(\bolds{\rho})$ in
(\ref{Hampluspert})
into two sums
\[
\sum_{k\leq\pi(\alpha(N-p+1))}\theta_k(\rho_{i_{1,k}},\ldots
,\rho_{i_{p,k}})
+
\sum_{k\leq\pi(\alpha(p-1))}\hat{\theta}_k(\rho_{j_{1,k}},\ldots
,\rho_{j_{p,k}}),
\]
where indices $i_{1,k},\ldots, i_{p,k}$ and $j_{1,k},\ldots, j_{p,k}$
are uniformly
distributed on $\{1,\ldots, N\}$. Therefore,
%
%e2.13 ###
\begin{equation}\label{this2}
\mathbb{E}\log\frac{Z_{N}}{Z_N'} =
\mathbb{E}\log\biggl\la\exp\sum_{k\leq\pi(\alpha(p-1))}\hat
{\theta}_k(\rho_{j_{1,k}}, \ldots,\rho_{j_{p,k}}) \biggr\ra'.
\end{equation}
Again Lemmas~\ref{LemP1} and~\ref{LemShift} imply that this converges to
\[
\mathbb{E}\log\mathbb{E}' \exp\sum_{k\leq\pi(\alpha(p-1))}
\hat{\theta}_k(s_{1,k},\ldots, s_{p,k}),
\]
and this finishes the proof of the lower bound.
\end{pf}

If we knew that $\mu\in\M$ is the unique limit of the sequence
$(\mu_N)$,
this would finish the proof of the first half of Theorem~\ref{ThFE}, since
$ \lim_{N\to\infty} F_N =
\lim_{N\to\infty} \mathbb{E}\log Z_{N+1}/Z_N
$
when the limit on the right exists. However, the proof of the general
case and
the second half of Theorem~\ref{ThFE} will require more work.
Before we move to the upper bound, let us record one more
consequence of the argument in Lemma~\ref{LemLB}. For $n\geq1$, let
us define
%
%e2.14 ###
\begin{equation}\label{Pn}\qquad
\P_n(\mu) = \log2 + \frac{1}{n} \mathbb{E}\log\mathbb{E}' \Av
_{\eps} \exp\sum_{i\leq n} A_i(\eps_i)
-
\frac{1}{n} \mathbb{E}\log\mathbb{E}' \exp\sum_{i\leq n} B_i.
\end{equation}
The following holds.
\begin{lemma} \label{LemInvn}
For all $\mu\in\M$, $\P_n(\mu)=\P(\mu)$ for all $n\geq1$.
\end{lemma}
\begin{pf}
We will only give a brief sketch since this will be
proved for all $\mu\in\M_{\mathrm{inv}}$
in Lemma~\ref{LemInvn2} below.
What we showed in the proof of Lemma~\ref{LemLB} is that if $\mu_N$
converges to $\mu$ over subsequence $(N_k)$, then $\mathbb{E}\log
Z_{N+1}/Z_N$ converges
to $\P(\mu)$ over the same subsequence. Similarly, one can show that,
given $n\geq1$,
over the same subsequence
\[
\frac{1}{n}( \mathbb{E}\log{Z_{N+n}} -\mathbb{E}\log{Z_{N}} ) \to
\P_n(\mu).\vadjust{\goodbreak}
\]
The only difference is that we split the terms in the Hamiltonian
$H_{N+n}(\vsi)$
into groups as in Lemma~\ref{LemShift2}, that is, instead of group (2)
we will have $n$ groups
each consisting of the terms with exactly one of the indices
$i_{1,k},\ldots, i_{p,k}$ equal to $N+l$ for $l=1,\ldots, n$.
On the other hand, if we write
\[
\frac{1}{n}( \mathbb{E}\log{Z_{N+n}} -\mathbb{E}\log{Z_{N}} )
=
\frac{1}{n}\sum_{l=1}^n( \mathbb{E}\log{Z_{N+l}} -\mathbb{E}\log
{Z_{N+l-1}} ),
\]
then repeating the proof of Lemma~\ref{LemLB} one can show that for
each term on
the right-hand side
\[
\lim_{N_k\to\infty} \mathbb{E}\log\frac{Z_{N_k+l}}{Z_{N_k+l-1}} =
\P(\mu),
\]
where instead of $\mu_{N_k}\to\mu$ one has to use that $\mu
_{N_{k}+l-1}\to\mu$
which holds by Lemma~\ref{LemShift2}. This finishes the proof.
\end{pf}

%s2.3 ###
\subsection{Upper bound and free energy}\label{SecUB}

Since the perturbation term in (\ref{Hampluspert}) does not affect the limit
of free energy, we will now ignore it and consider free energy $F_N$ defined
for the original unperturbed Hamiltonian (\ref{Ham}). Recall $A_i(\eps
)$ and $B_i$
defined in (\ref{Ai}) and (\ref{Bi}).
\begin{lemma}\label{LemUpper}
For any function $\sigma\dvtx[0,1]^4\to\{-1,+1\}$ we have
%
%e2.15 ###
\begin{equation}\label{upperbound}\qquad
F_N \leq \log2 + \frac{1}{N} \mathbb{E}\log\mathbb{E}' \Av_{\eps} \exp
\sum_{i\leq N} A_i(\eps_i) - \frac{1}{N} \mathbb{E}\log\mathbb{E}'
\exp\sum_{i\leq N} B_i.
\end{equation}
\end{lemma}
\begin{Remark*}
In general, this upper bound does not decouple and
depends on $N$ since
all $s_{i,k,l}$ and $\hat{s}_{i,k,l}$ defined in (\ref{s}) and (\ref
{hats}) depend on
the same variable $u$ in the second coordinate. We will see that the
proof of the upper bound
(\ref{upperbound}) does not to work if one tries to replace $u$ by
independent copies $u_i$ in
the definition of $A_i(\eps)$ and $B_i$. For $\sigma=\sigma_\mu$
for $\mu\in\M$,
Lemma~\ref{LemInvn} implies that this upper bound does not depend on
$N$ and, thus,
$F_N\leq\P(\mu)$. Together with the lower bound of Lemma \ref
{LemLB} this proves that
\[
\lim_{N\to\infty} F_N = \inf_{\mu\in\M} \P(\mu).
\]
To prove the second part of Theorem~\ref{ThFE}, we will show in Lemma
\ref{LemInvn2} below
that the invariance properties in (\ref{SC}) imply that $\P_n(\mu) =
\P(\mu)$ for $\mu\in\M_{\mathrm{inv}}$
as well which will finish the proof of Theorem~\ref{ThFE}.
\end{Remark*}
\begin{pf*}{Proof of Lemma~\ref{LemUpper}}
A proof by interpolation is a slight modification of the proof in~\cite{PT}.
For $t\in[0,1]$, let us define similarly to (\ref{Ai}) and (\ref{Bi})
%
%e2.16 ###
\begin{equation}\label{Ait}
A_i^t(\eps)=\sum_{k\leq\pi_i((1-t)p\alpha)} \theta_{k,i}(\eps,
s_{i,k,1},\ldots,s_{i,k,p-1})
\end{equation}
and
%
%e2.17 ###
\begin{equation}\label{Bit}
B_i^t = \sum_{k\leq\pi_i(t(p-1)\alpha)} \hat{\theta}_{k,i}(\hat
{s}_{i,k,1},\ldots,\hat{s}_{i,k,p}).
\end{equation}
Consider an interpolating Hamiltonian
%
%e2.18 ###
\begin{equation}\label{inter-rs}
-H_{N,t}(\vsi)=\sum_{k\leq\pi(t\alpha N)}\theta_k(\sigma
_{i_{1,k}},\ldots,\sigma_{i_{p,k}}) + \sum_{i\leq
N}A_i^t(\sigma_i)+\sum_{i\leq N} B_i^t
\end{equation}
and let
\[
\varphi(t)=\frac{1}{N} \mathbb{E}\log\mathbb{E}'\sum_{\vsi\in
\Sigma_N} \exp(-H_{N,t}(\vsi)).
\]
Since, clearly,
\[
\varphi(1) = F_N + \frac{1}{N} \mathbb{E}\log\mathbb{E}' \exp
\sum_{i\leq N} B_i
\]
and
\[
\varphi(0) = \log2 + \frac{1}{N} \mathbb{E}\log\mathbb{E}' \Av
_{\eps} \exp\sum_{i\leq N} A_i(\eps_i),
\]
it remains to prove that $\varphi'(t)\leq0$.
Let us consider the partition function
\[
Z=\sum_{\vsi\in\Sigma_N} \exp(-H_{N,t}(\vsi))
\]
and define
\[
Z_{m}=Z|_{\pi(t\alpha N) = m},\qquad
Z_{i,m}^A=Z|_{\pi_i((1-t)p\alpha) = m}
\quad\mbox{and}\quad
Z_{i,m}^B=Z|_{\pi_i(t (p-1)\alpha) = m}.
\]
If we denote the Poisson p.f. as
$\pi(\lambda,k)=(\lambda^k/k!) e^{-\lambda}$,
then
\[
\mathbb{E}\log\mathbb{E}' Z = \sum_{m\geq0} \pi(t \alpha N,m)
\mathbb{E}\log\mathbb{E}' Z_m
\]
and, for any $i \leq N$,
\[
\mathbb{E}\log\mathbb{E}' Z = \sum_{m\geq0} \pi\bigl((1-t)p\alpha, m\bigr)
\mathbb{E}\log\mathbb{E}' Z_{i,m}^A
\]
and
\[
\mathbb{E}\log\mathbb{E}' Z = \sum_{m\geq0} \pi\bigl(t(p-1)\alpha, m\bigr)
\mathbb{E}\log\mathbb{E}' Z_{i,m}^B.
\]
Therefore, we can write
%
%e2.19 ###
\begin{eqnarray}\label{dersteps}\quad
\varphi'(t)
&=&
\sum_{m\geq0} \frac{\partial\pi(t\alpha N,m)}{\partial t}
\frac{1}{N}  \mathbb{E}\log\mathbb{E}' Z_m\nonumber\\[-1pt]
&&{}+
\sum_{i\leq N} \sum_{m\geq0}
\frac{\partial\pi((1-t)p\alpha, m)}{\partial t}
\frac{1}{N}  \mathbb{E}\log\mathbb{E}' Z_{i,m}^A
\nonumber\\[-1pt]
&&{}
+
\sum_{i\leq N} \sum_{m\geq0}
\frac{\partial\pi(t(p-1)\alpha, m)}{\partial t}
\frac{1}{N}  \mathbb{E}\log\mathbb{E}' Z_{i,m}^B
\nonumber\\[-1pt]
&=&
\alpha\sum_{m\geq0}
\bigl(\pi(t\alpha N,m-1)I(m\geq1) - \pi(t\alpha N,m) \bigr)
\mathbb{E}\log\mathbb{E}' Z_m
\nonumber\\[-1pt]
&&{}
- p \alpha
\frac{1}{N}
\sum_{i\leq N} \sum_{m\geq0}
\bigl(\pi\bigl((1-t)p\alpha, m-1\bigr)I(m\geq1)\nonumber\\[-1pt]
&&\hspace*{127pt}{} - \pi\bigl((1-t)p\alpha,m\bigr)\bigr)
\mathbb{E}\log\mathbb{E}' Z_{i,m}^A
\nonumber\\[-9pt]\\[-9pt]
&&{}
+ (p-1) \alpha
\frac{1}{N}
\sum_{i\leq N} \sum_{m\geq0}
\bigl(\pi\bigl(t (p-1)\alpha, m-1\bigr)I(m\geq1)\nonumber\\[-1pt]
&&\hspace*{154pt}{} -
\pi\bigl(t(p-1)\alpha,m\bigr)\bigr)
\mathbb{E}\log\mathbb{E}' Z_{i,m}^B
\nonumber\\[-1pt]
&=&
\alpha\sum_{m\geq0}
\pi(t\alpha N,m) \mathbb{E}\log(\mathbb{E}' Z_{m+1}/\mathbb{E}' Z_m)
\nonumber\\[-1pt]
&&{}
- p\alpha
\frac{1}{N}
\sum_{i\leq N} \sum_{m\geq0}
\pi\bigl((1-t)p \alpha, m\bigr) \mathbb{E}\log(\mathbb{E}'
Z_{i,m+1}^A/\mathbb{E}' Z_{i,m}^A)
\nonumber\\[-1pt]
&&{}
+ (p-1)\alpha
\frac{1}{N}
\sum_{i\leq N} \sum_{m\geq0}
\pi\bigl(t(p-1) \alpha, m\bigr) \mathbb{E}\log(\mathbb{E}'
Z_{i,m+1}^B/\mathbb{E}' Z_{i,m}^B)
\nonumber\\[-1pt]
&=&
\alpha\mathbb{E}\log\frac{\mathbb{E}' Z_{+1}}{\mathbb{E}' Z}
- p\alpha
\frac{1}{N}
\sum_{i\leq N} \mathbb{E}\log\frac{\mathbb{E}' Z_{i,+1}^A}{\mathbb
{E}' Z}
+ (p-1)\alpha
\mathbb{E}\log\frac{\mathbb{E}' Z_{+1}^B}{\mathbb{E}' Z},\nonumber
\end{eqnarray}
where $Z_{+1}$, $Z_{i,+1}^A$ and $Z_{+1}^B$ contain one extra term in
the Hamiltonian in the corresponding
Poisson sum. Namely,
\begin{eqnarray*}
Z_{+1} &=& \sum_{\vsi\in\Sigma_N} \exp\theta(\sigma_{i_{1}},\ldots
,\sigma_{i_{p}}) \exp(-H_{N,t}(\vsi)),
\\[-1pt]
Z_{i,+1}^A &=& \sum_{\vsi\in\Sigma_N} \exp\theta(\sigma
_{i},s_1,\ldots, s_{p-1}) \exp(-H_{N,t}(\vsi)),
\\[-1pt]
Z_{+1}^B &=& \sum_{\vsi\in\Sigma_N} \exp\theta(s_1,\ldots, s_{p})
\exp(-H_{N,t}(\vsi)),
\end{eqnarray*}
where random function $\theta$ and indices $i_1,\ldots, i_p$ uniform
on $\{1,\ldots, N\}$
are independent of the randomness of the Hamiltonian $H_{N,t}$.
If, for a function $f$ of $\vsi, u$ and $(x)$, we denote by $\la f\ra_t$
the Gibbs average
\[
\la f\ra_t = \frac{1}{\mathbb{E}' Z}  \mathbb{E}' \sum_{\vsi\in
\Sigma_N} f \exp(-H_{N,t}(\vsi)),
\]
then (\ref{dersteps}) can be rewritten as
%
%e2.20 ###
\begin{eqnarray}\label{neg}
&&
\alpha  \mathbb{E}\log
\langle\exp\theta(\sigma_{i_{1}},\ldots,\sigma
_{i_{p}})\rangle_t\nonumber\\
&&\qquad{}-p\alpha \frac{1}{N}\sum_{i\leq N}
\mathbb{E}\log
\langle\exp\theta(\sigma_{i},s_1,\ldots, s_{p-1})
\rangle_t
\\
&&\qquad{}
+ (p-1) \alpha
\mathbb{E}\log
\langle\exp{\theta}({s}_1,\ldots, {s}_{p})\rangle_t.
\nonumber
\end{eqnarray}
By assumptions (\ref{condition1}) and (\ref{condition3}) we can write
\begin{eqnarray*}
\log\langle\exp\theta(\sigma_{i_{1}},\ldots,\sigma
_{i_{p}})\rangle_t
&=&
\log a +
\log\bigl(1+b\langle f_1(\sigma_{i_1})\cdots f_p(\sigma
_{i_p})\rangle_t\bigr)
\\
&=& \log a -
\sum_{n\geq1} \frac{(-b)^n}{n}\langle f_1(\sigma_{i_1})\cdots
f_p(\sigma_{i_p}) \rangle_t^n.
\end{eqnarray*}
Using replicas $\vsi^l, u_l$ and $(x^l)$, we can write
\[
\langle f_1(\sigma_{i_1})\cdots f_p(\sigma_{i_p}) \rangle
_t^n =
\biggl\langle\prod_{l\leq n} f_1(\sigma_{i_1}^l)\cdots f_p(\sigma
_{i_p}^l) \biggr\rangle_t
\]
and thus
\[
\frac{1}{N^p}\sum_{i_1,\ldots,i_p\leq N} \langle f_1(\sigma
_{i_1})\cdots
f_p(\sigma_{i_p}) \rangle_t^n =
\biggl\langle\prod_{j\leq p} A_{j,n} \biggr\rangle_t,
\]
where
\[
A_{j,n}=A_{j,n}(\vsi^1,\ldots, \vsi^n) = \frac{1}{N}\sum_{i\leq N}
\prod_{l\leq n}f_j(\sigma_i^l).
\]
Denote by $\mathbb{E}_0$ the expectation in $f_1,\ldots,f_p$.
Since $f_1,\ldots,f_p$ are i.i.d. and independent of the randomness in
$\la\cdot\ra_t$,
\[
\mathbb{E}_0 \biggl\la\prod_{j\leq p} A_{j,n}\biggr\ra_t
= \biggl\la\mathbb{E}_0 \prod_{j\leq p} A_{j,n}\biggr\ra_t
=\la B_n^p\ra_t,
\]
where $B_n = \mathbb{E}_0 A_{j,n}$.
Therefore, since we also assumed that $b$ is independent of
$f_1,\ldots,f_p$,
%
%e2.21 ###
\begin{equation}\label{neg1}
\mathbb{E}_0 \frac{1}{N^p}
\sum_{i_1,\ldots,i_p\leq N}
\log\langle\exp\theta(\sigma_{i_{1}},\ldots,\sigma
_{i_{p}})\rangle_t
= \mathbb{E}_0 \log a
- \sum_{n\geq1} \frac{(-b)^n}{n} \langle B_n^p \rangle_t.\hspace*{-32pt}
\end{equation}
A similar analysis applies to the second term in (\ref{neg}),
\begin{eqnarray*}
&&\log\langle\exp\theta(\sigma_{i},s_1,\ldots, s_{p-1})
\rangle_t\\
&&\qquad= \log a -\sum_{n\geq1}\frac{(-b)^n}{n}
\biggl\langle f_p(\sigma_i) \prod_{j\leq p-1} f_j(s_j) \biggr\rangle_t^n
\\
&&\qquad=
\log a
-\sum_{n\geq1}\frac{(-b)^n}{n}
\biggl\langle\prod_{l\leq n} f_p(\sigma_i^l)
\prod_{l\leq n} \prod_{j\leq p-1}
f_j(s_j^l)
\biggr\rangle_t,
\end{eqnarray*}
where in the last equality we again used replicas $\vsi^l, u_l$ and $(x^l)$;
for example, compared to (\ref{s}), $s_j^l$ is now defined by $s^l_{j}
= \sigma(w,u_l, v_{j},x^l_{j})$.
Thus,
\[
\frac{1}{N} \sum_{i\leq N}\log\langle\exp\theta(\sigma
_{i},s_1,\ldots, s_{p-1}) \rangle_t
=
\log a
-\sum_{n\geq1}\frac{(-b)^n}{n}
\biggl\langle A_{p,n}
\prod_{j\leq p-1} \prod_{l\leq n}
f_j(s_j^l)
\biggr\rangle_t.
\]
[\textit{Note}: It was crucial here that $s_j^l$ do not depend on $i$
through independent
copies $u_i$ rather than the same $u$. It is tempting to define the
interpolation (\ref{inter-rs})
by using independent $u_i$ for $i\leq N$ since this would make the
upper bound in (\ref{upperbound}) decouple, but the proof would break
down at this step.]
In addition to $f_1,\ldots, f_p$, let $\mathbb{E}_0$ also denote the
expectation in $(v_j)$
and $(x_j^l)$ in $s_j^l$, but not in sequences $(v), (x)$ in the
randomness of $\la\cdot\ra_t$.
Then,
%
%e2.22 ###
\begin{eqnarray}\label{neg2}
&&\mathbb{E}_0 \frac{1}{N} \sum_{i\leq N}\log\langle\exp
\theta(\sigma_{i},s_1,\ldots, s_{p-1}) \rangle_t\nonumber\\[-8pt]\\[-8pt]
&&\qquad = \mathbb{E}_0 \log a
-\sum_{n\geq1}\frac{(-b)^n}{n} \langle B_n
(C_n)^{p-1}\rangle_t,\nonumber
\end{eqnarray}
where
\[
C_n = C_n(w,u_1,\ldots, u_n)
= \mathbb{E}_0 \prod_{l\leq n} f_j(s_j^l)
=
\mathbb{E}_0 \prod_{l\leq n} f_j(\sigma(w,u_l,v_j,x_j^l))
\]
obviously does not depend on $j$. Finally, in an absolutely similar manner
%
%e2.23 ###
\begin{equation}\label{neg3}
\mathbb{E}_0 \log\langle\exp{\theta}({s}_1,\ldots, {s}_{p})\rangle_t =
\mathbb{E}_0 \log a -\sum_{n\geq1}\frac{(-b)^n}{n} \langle (C_n)^p
\rangle_t.
\end{equation}
Combining (\ref{neg1}), (\ref{neg2}) and (\ref{neg3}) we see that
(\ref{neg}) can be written as
%
%e2.24 ###
\begin{equation}\label{negdone}
-\alpha\sum_{n\geq1}\frac{\mathbb{E}(-b)^n}{n} \mathbb{E}\langle B_n^p
-p B_n C_n^{p-1} + (p-1) (C_n)^p \rangle_t \leq0,
\end{equation}
which holds true using condition (\ref{condition2}) and the fact that
$x^p - pxy^{p-1} +(p-1)y^p\geq0$ for all $x,y\in\Reals$ for even
$p\geq2$.
This finishes the proof of the upper bound.
\end{pf*}

Before proving the invariance properties of Theorem~\ref{ThSC} let us
finish the proof
of Theorem~\ref{ThFE} by showing that for invariant measures $\M
_{\mathrm{inv}}$ the upper bound decouples.
\begin{lemma} \label{LemInvn2}
For all $\mu\in\M_{\mathrm{inv}}$, $\P_n(\mu)=\P(\mu)$ for all $n\geq1$.
\end{lemma}
\begin{pf}
If we recall $A_i$ defined in (\ref{Aiav}), then we can rewrite (\ref
{Pn}) as
%
%e2.25 ###
\begin{equation}\label{Pn2}
\P_n(\mu) = \log2 + \frac{1}{n} \mathbb{E}\log\mathbb{E}' \exp
\sum_{i\leq n} A_i
-\frac{1}{n} \mathbb{E}\log\mathbb{E}' \exp\sum_{i\leq n} B_i.
\end{equation}
The result will follow if we show that for any $n\geq1$,
%
%e2.26 ###
\begin{equation}\label{Arec}
\mathbb{E}\log\frac{\mathbb{E}' \exp\sum_{i\leq n+1} A_i}{\mathbb {E}'
\exp\sum_{i\leq n} A_i} = \mathbb{E}\log\mathbb{E}' \exp A_{n+1}
\end{equation}
and
%
%e2.27 ###
\begin{equation}\label{Brec}
\mathbb{E}\log\frac{\mathbb{E}' \exp\sum_{i\leq n+1} B_i}{\mathbb
{E}' \exp\sum_{i\leq n} B_i}
=\mathbb{E}\log\mathbb{E}' \exp B_{n+1}.
\end{equation}
To prove this we will use the invariance properties (\ref{ASC}) and
(\ref{preBSC}).
If in (\ref{preBSC}) we choose $r$ to be a Poisson r.v. with mean $n
(p-1) \alpha$, then it becomes
%
%e2.28 ###
\begin{equation}\label{BSC}
\mathbb{E}\prod_{l\leq q}\mathbb{E}'\prod_{i\in C_l} s_i
=
\mathbb{E}\frac{\prod_{l\leq q} \mathbb{E}' \prod_{i\in C_l} s_i
\exp\sum_{i\leq n} B_i}
{(\mathbb{E}' \exp\sum_{i\leq n} B_i)^q}.
\end{equation}
We will only show how (\ref{ASC}) implies (\ref{Arec}) since the
proof that
(\ref{BSC}) implies (\ref{Brec}) is exactly the same. We only need to prove
(\ref{Arec}) conditionally on the Poisson r.v. $\pi_{n+1}(p\alpha)$
and functions
$(\theta_{k,n+1})$ in the definition of $A_{n+1}$,
%
%e2.29 ###
\begin{equation}\label{Airemind}
\exp A_{n+1} = \Av_\eps\exp\sum_{k\leq\pi_{n+1}(p\alpha)}
\theta_{k}(\eps, s_{1,n+1,k},\ldots,s_{p-1,n+1,k}),
\end{equation}
since we can control these functions
uniformly with high probability using condition (\ref{condition4}).
Approximating the
logarithm by polynomials, in order to prove (\ref{Arec}), it is enough
to prove that
%
%e2.30 ###
\begin{equation}\label{ASCq}
\mathbb{E}\biggl(\frac{\mathbb{E}' \exp A_{n+1} \exp\sum_{i\leq n}
A_i}{\mathbb{E}' \exp\sum_{i\leq n} A_i}\biggr)^q
=\mathbb{E}(\mathbb{E}' \exp A_{n+1})^q
\end{equation}
for all $q\geq1$. Condition (\ref{condition1}) implies that the
right-hand side of (\ref{Airemind})
is a polynomial of spins $(s_{j,n+1,k})$ for $k\leq\pi_{n+1}(p\alpha
)$ and $j\leq p-1$,
and, therefore, (\ref{ASCq}) is obviously implied by (\ref{ASC}) if
we simply enumerate spins
$(s_{j,n+1,k})$ as spins $(s_i)$ for $n+1\leq i\leq m$ by choosing $m$
large enough.
Averaging over random $\pi_{n+1}(p\alpha)$ and $(\theta_{k,n+1})$
proves (\ref{ASCq}) and finishes the proof.
\end{pf}

Let us note that, similarly, (\ref{preBSC}) implies
\[
\mathbb{E}\log\frac{\mathbb{E}' \exp\sum_{i\leq n+1} \hat{\theta
}_i(\hat{s}_{1,i},\ldots,\hat{s}_{p,i})}
{\mathbb{E}' \exp\sum_{i\leq n} \hat{\theta}_i(\hat
{s}_{1,i},\ldots,\hat{s}_{p,i})}
=
\mathbb{E}\log\mathbb{E}' \exp\hat{\theta}_1(\hat
{s}_{1,1},\ldots,\hat{s}_{p,1}),
\]
which obviously implies (\ref{Plast}), that is,
\begin{eqnarray*}
\mathbb{E}\log\mathbb{E}' \exp B
&=&
\mathbb{E}\log\mathbb{E}' \exp\sum_{k\leq\pi((p-1) \alpha)}
\theta_{k}(s_{1,k},\ldots,s_{p,k})
\\
&=&
(p-1)\alpha \mathbb{E}\log\mathbb{E}' \exp\theta(s_1,\ldots,s_p).
\end{eqnarray*}

%s2.4 ###
\subsection{Invariance and self-consistency equations}\label{SecSC}

\mbox{}

\begin{pf*}{Proof of Theorem~\ref{ThSC}}
Let $h = \prod_{l\leq q} h_l$ where $h_l = \prod_{j\in
C_l}\sigma_j^l$.\vspace*{1pt}
Consider $\mu\in\M$ which is a limit
of $\mu_N$ over some subsequence $(N_k)$.
Using Lemma~\ref{LemShift2}, the left-hand side of (\ref{SC}) is the
limit of
$\mathbb{E}\la h \ra_{N+n}$ over subsequence $(N_k)$.
The right-hand side of (\ref{SC}) will appear as a similar limit once
we rewrite this joint moment
of spins using cavity coordinates and ``borrowing'' some terms in the
Gibbs measure from
the Hamiltonian (\ref{Hampluspert}). The spins with coordinates $i\leq
n$ will play the role
of cavity coordinates. Let us separate the $\pi(\alpha(N+n))$ terms
in the first sum
%
%e2.31 ###
\begin{equation}\label{firstsumI}
\sum_{k\leq\pi(\alpha(N+n))}\theta_k(\sigma_{i_{1,k}}, \ldots
,\sigma_{i_{p,k}})
\end{equation}
in (\ref{Hampluspert}) in the Hamiltonian $H_{N+n}$ into three groups:

\begin{longlist}[($2j$)]
\item[(1)]
terms for $k$ such that all indices $i_{1,k},\ldots, i_{p,k}>
n$;

For $1\leq j\leq n$:

\item[($2j$)] terms with exactly one of indices $i_{1,k},\ldots, i_{p,k}$
equal to $j$ and all others~$>$~$n$;

\item[(3)] terms with at least two of indices $i_{1,k},\ldots, i_{p,k}\leq
n$.
\end{longlist}

The probabilities that a term is of these three type are
\[
p_1=\biggl(\frac{N}{N+n} \biggr)^p,\qquad
p_{2,j} = p\frac{1}{N+n} \biggl(\frac{N}{N+n}\biggr)^{p-1},\qquad
p_3=1-p_1-\sum_{l\leq n}p_{2,l}.
\]
Therefore, the number of terms in these groups are independent
Poisson random variables with means
\begin{eqnarray*}
\alpha(N+n) p_1 &=& \alpha(N +n- np) + O(N^{-1}),\\
\alpha(N+n) p_{2,j} &=& \alpha p + O(N^{-1}),\\
\alpha(N+n) p_3 &=& O(N^{-1}).
\end{eqnarray*}
We can redefine the number of terms in each group to be exactly of
means $\alpha(N+n-np),
\alpha p$ and $0$ since asymptotically it does not affect $\mathbb
{E}\la h\ra_{N+n}$.
Thus, if we write $\vsi=(\bolds{\eps},\bolds{\rho})\in\Sigma_{N+n}$
for the first
the first $n$ cavity\vadjust{\goodbreak} coordinates $\bolds{\eps}=(\eps_1,\ldots, \eps
_{n})$ and the last
$N$ coordinates $\bolds{\rho}=(\rho_1,\ldots,\rho_{N})$, then (\ref{firstsumI})
can be replaced with
%
%e2.32 ###
\begin{eqnarray}\label{firsttermI}
&&\sum_{k\leq\pi(\alpha(N+n-np))}\theta_k(\rho_{i_{1,k}},
\ldots,\rho_{i_{p,k}})\nonumber\\[-8pt]\\[-8pt]
&&\qquad{}+\sum_{j\leq n}
\sum_{k\leq\pi_j(\alpha p)}\theta_{k,j}(\eps_j, \rho_{i_{1,k,j}},
\ldots,\rho_{i_{p-1,k,j}}),\nonumber
\end{eqnarray}
where indices $i_{1,k},\ldots, i_{p,k}$ and $i_{1,k,j},\ldots,
i_{p-1,k,j}$ are all
uniformly distributed on $\{1,\ldots, N\}$. Let us now consider the
perturbation term
in (\ref{Hampluspert}),
%
%e2.33 ###
\begin{equation}\label{pertnI}
\sum_{l\leq\pi(c_{N+n})} \log\Av_\eps \exp\sum_{k\leq\hat{\pi}_l(\alpha
p)} \hat{\theta}_{k,l}(\eps, \sigma_{j_{1,k,l}},
\ldots,\sigma_{j_{p-1,k,l}}) ,
\end{equation}
where $j_{1,k,l},\ldots, j_{p-1,k,l}$ are uniformly distributed on $\{
1,\ldots, N+n\}$.
Here, we used independent copies $\hat{\pi}_l$ and $\hat{\theta
}_{k,l}$ since
$\pi_j$ and $\theta_{k,j}$ were already used in (\ref{firsttermI}).
The expected number of these indices that belong to $\{1,\ldots, n\}$ is
$c_{N+n}\alpha p(p-1)n/N\to0$ which means that with high probability
all indices belong to $\{n+1,\ldots, N+n\}$. As a result,
asymptotically $\mathbb{E}\la h\ra_{N+n}$
will not be affected if we replace the perturbation term (\ref
{pertnI}) with
%
%e2.34 ###
\begin{equation}\label{pertn2I}
\sum_{l\leq\pi(c_{N+n})} \log\Av_\eps \exp\sum_{k\leq\hat{\pi}_l(\alpha
p)} \hat{\theta}_{k,l}(\eps, \rho_{j_{1,k,l}},
\ldots,\rho_{j_{p-1,k,l}}),
\end{equation}
where $j_{1,k,l},\ldots, j_{p-1,k,l}$ are uniformly distributed on $\{
1,\ldots, N\}$.
Thus, we can assume from now on that $\mathbb{E}\la h\ra_{N+n}$ is
computed with respect to the Hamiltonian
which is the sum of (\ref{firsttermI}) and (\ref{pertn2I}). If $\la
\cdot\ra_N'$ denotes the Gibbs average
on $\Sigma_N$ with respect to the Hamiltonian
%
%e2.35 ###
\begin{eqnarray}\label{HprimeI}\qquad
-H_N'(\bolds{\rho})
&=&
\sum_{k\leq\pi(\alpha(N+n-np))}\theta_k(\rho_{i_{1,k}},
\ldots,\rho_{i_{p,k}})
\nonumber\\[-8pt]\\[-8pt]
&&{}
+
\sum_{l\leq\pi(c_{N+n})}
\log\Av_\eps
\exp\sum_{k\leq\hat{\pi}_l(\alpha p)}
\hat{\theta}_{k,l}(\eps, \rho_{j_{1,k,l}}, \ldots,\rho_{j_{p-1,k,l}}),
\nonumber
\end{eqnarray}
then we can write
%
%e2.36 ###
\begin{equation}\label{reprI}
\mathbb{E}\la h\ra_{N+n} = \mathbb{E} \frac{\prod_{l\leq q}
U_{N,l}}{V_{N}^q},
\end{equation}
where
\[
U_{N,l} = \biggl\la \Av_\eps  h_l(\bolds{\eps}, \bolds{\rho}) \exp\sum_{j\leq n}
\sum_{k\leq\pi_j(\alpha p)}\theta_{k,j}(\eps_j, \rho_{i_{1,k,j}},
\ldots,\rho_{i_{p-1,k,j}})\biggr\ra_N'
\]
and
\[
V_{N} = \biggl\la \Av_\eps\exp\sum_{j\leq n} \sum_{k\leq\pi_j(\alpha
p)}\theta_{k,j}(\eps_j, \rho_{i_{1,k,j}},
\ldots,\rho_{i_{p-1,k,j}})\biggr\ra_N'.\vadjust{\goodbreak}
\]
Finally, given $r\geq1$, let us borrow $r$ terms from the first sum in
(\ref{HprimeI})
by splitting the last $r$ terms and replacing the first sum in (\ref
{HprimeI}) with
\[
\sum_{k\leq\pi(\alpha(N+n-np))-r}\theta_k(\rho_{i_{1,k}},
\ldots,\rho_{i_{p,k}}) + \sum_{k\leq r}\hat{\theta}_k(\rho_{j_{1,k}},
\ldots,\rho_{j_{p,k}}).
\]
Here we ignore the negligible event when $\pi(\alpha(N+n-np))<r$.
If we define
%
%e2.37 ###
\begin{eqnarray}\label{HprimeI2}\qquad
-H_N''(\bolds{\rho})
&=&
\sum_{k\leq\pi(\alpha(N+n-np))-r}\theta_k(\rho_{i_{1,k}},
\ldots,\rho_{i_{p,k}})
\nonumber\\[-8pt]\\[-8pt]
&&{}
+
\sum_{l\leq\pi(c_{N+n})}
\log\Av_\eps
\exp\sum_{k\leq\hat{\pi}_l(\alpha p)}
\hat{\theta}_{k,l}(\eps, \rho_{j_{1,k,l}}, \ldots,\rho_{j_{p-1,k,l}})
\nonumber
\end{eqnarray}
and let $\la\cdot\ra_N''$ denote the Gibbs average on $\Sigma_N$
with respect to this
Hamiltonian then $U_{N,l}/V_{N} = U_{N,l}'/V_{N}'$ where
\begin{eqnarray*}
U_{N,l}' &=& \biggl\la
\Av_\eps  h_l(\bolds{\eps}, \bolds{\rho}) \exp\sum_{j\leq n}
\sum_{k\leq\pi_j(\alpha p)}\theta_{k,j}(\eps_j, \rho_{i_{1,k,j}},
\ldots,\rho_{i_{p-1,k,j}})\\
&&\qquad\hspace*{103.5pt}{}\times\exp\sum_{k\leq r}\hat{\theta}_k(\rho_{j_{1,k}},
\ldots,\rho_{j_{p,k}})
\biggr\ra_N''
\end{eqnarray*}
and
\begin{eqnarray*}
V_{N}'' &=& \biggl\la
\Av_\eps\exp\sum_{j\leq n}
\sum_{k\leq\pi_j(\alpha p)}\theta_{k,j}(\eps_j, \rho_{i_{1,k,j}},
\ldots,\rho_{i_{p-1,k,j}})\\
&&\hspace*{89pt}{}\times\exp\sum_{k\leq r}\hat{\theta}_k(\rho_{j_{1,k}},
\ldots,\rho_{j_{p,k}})
\biggr\ra_N''.
\end{eqnarray*}
By Lemma~\ref{LemShift}, the distribution of spins under the annealed
Gibbs measure
$\mathbb{E}\la\cdot\ra_N''$ corresponding to the Hamiltonian
$H_N''(\bolds{\rho})$ still converges
to $\mu$ over the subsequence~$(N_k)$. Conditionally on $(\pi
_j(\alpha p))$,
$(\theta_{k,j})$, $(\hat{\theta}_k)$ and on the event that all indices
$i_{1,k,j}, \ldots,i_{p-1,k,j}$ and $j_{1,k},\ldots,j_{p,k} $ are
different, Lemma~\ref{LemP1}
implies that the right-hand side of (\ref{reprI}) converges over
subsequence\vspace*{2pt} $(N_k)$
to $\mathbb{E}\prod_{l\leq q} U_{l}/V^q$ where $(U_l)$ and $V$ are
defined in (\ref{Ul}) and (\ref{Vl})
only now conditionally on the above sequences.
Since asymptotically all indices are different with high probability,
the same convergence holds
unconditionally, and this completes the proof.
\end{pf*}

%s3 ###
\section{Sherrington--Kirkpatrick model}\label{SecSKK}

%s3.1 ###
\subsection{Properties of convergence}\label{SecPCSK}

Of course, Lemma~\ref{LemP1} still holds since it does not really
depend on the model.
However, the role of this lemma in the\vadjust{\goodbreak} Sherrington--Kirkpatrick model
will be played by the statement that
we made at the beginning of the introduction which we now record for
the reference.
\begin{lemma}\label{LemP1SK}
The joint distribution of spins $(\sigma_i^l)$ and multi-overlaps
(\ref{multioverlapN}) converges to
the joint distribution of spins (\ref{sigma}) and multi-overlaps (\ref
{multioverlap}) over any subsequence
along which $\mu_N$ converges to $\mu$.
\end{lemma}

Lemma~\ref{LemShift} also has a straightforward analog for the
Sherrington--Kirkpatrick model.
Let $\la\cdot\ra$ denote the Gibbs average with respect to the sum
of an arbitrary Hamiltonian
on $\Sigma_N$ and a perturbation term (\ref{HampertSK}), and let $\la
\cdot\ra'$ denote
the Gibbs average corresponding to the sum of the same arbitrary Hamiltonian
and a perturbation as in (\ref{HampertSK}), only with the number of
terms replaced by
$\pi(c_N) + n$ instead of $\pi(c_N)$ in the first sum and $\pi
'(c_N)+m$ instead of $\pi'(c_N)$
in the second sum, for any finite $m,n\geq1$. Then the following holds.

\begin{lemma}\label{LemShiftSK}
For any bounded function $h$ of finitely many spins, or finitely many
multi-overlaps, we have
%
%e3.1 ###
\begin{equation}\label{ShiftSK}
{\lim_{N\to\infty}} | \mathbb{E}\la h \ra' -\mathbb{E}\la h\ra|=0.
\end{equation}
\end{lemma}

The proof is exactly the same as in Lemma~\ref{LemShift}.
The role of the perturbation (\ref{HampertSK}) will finally start
becoming clear in the following
exact analog of Lemma~\ref{LemShift2}.
\begin{lemma}\label{LemShift2SK}
If $\mu_N$ converges to $\mu$ over subsequence $(N_k)$, then it also converges
to $\mu$ over subsequence $(N_k+n)$ for any $n\geq1$.
\end{lemma}
\begin{pf}
We will\vspace*{1pt} show that the joint moments of spins converge
to the same limit over
subsequences that differ by a finite shift $n$. Let $h=\prod_{j\leq q}
h_j$ where
$h_j = \prod_{i\in C_j} \sigma_i^j$ over some\vspace*{1pt} finite sets of spin
coordinates $C_j$.
Let us denote by $\la\cdot\ra_N$ the
Gibbs average with respect to the Hamiltonian (\ref{Hampluspert})
defined on $N$ coordinates.
We will show that
\[
{\lim_{N\to\infty}} |\mathbb{E}\la h\ra_{N+n} - \mathbb{E}\la h\ra
_N| = 0.
\]
Let us rewrite $\mathbb{E}\la h\ra_{N+n}$ by treating the last $n$
coordinates as cavity coordinates.
Let us write $\vsi=(\bolds{\rho},\bolds{\eps})\in\Sigma_{N+n}$ for
the first
$N$ coordinates $\bolds{\rho}=(\rho_1,\ldots,\rho_{N})$ and the last $n$
cavity coordinates $\bolds{\eps}=(\eps_1,\ldots, \eps_{n})$ and
rewrite (\ref{HamSK}) as
%
%e3.2 ###
\begin{equation}\label{firstsumSK}
-H_{N+n}(\bolds{\rho}) + \sum_{i\leq n}\eps_i Z_i(\bolds{\rho}) +
\delta(\vsi),
\end{equation}
where we define (slightly abusing notations)
%
%e3.3 ###
\begin{equation}\label{HamSKcav}
-H_{N+n}(\bolds{\rho}) :=\sum_{p\geq1}\frac{\beta_p}{(N+n)^{(p-1)/2}}
\sum_{1\leq i_1,\ldots,i_p\leq N}g_{i_1,\ldots,i_p} \rho
_{i_1}\cdots\rho_{i_p};\vadjust{\goodbreak}
\end{equation}
the term $\eps_i Z_i(\bolds{\rho})$ consists of all terms in (\ref
{HamSK}) with only one factor
$\eps_i$ from $\bolds{\eps}$ present, and the last term $\delta$ is
the sum of terms with at least
two factors in $\bolds{\eps}$. It is easy to check that
\[
\mathbb{E}Z_i(\bolds{\rho}^{1}) Z_i(\bolds{\rho}^{2})
= \xi'(R(\bolds{\rho}^{1},\bolds{\rho}^{2})) + o_N(1)
\]
uniformly over all $\vrho^1, \vrho^2$, and the covariance of $\delta
(\vsi)$ is also
of small order uniformly over $\vsi^1,\vsi^2$. By the usual Gaussian
interpolation
one can therefore redefine the Hamiltonian $H_{N+n}(\vsi)$ by
%
%e3.4 ###
\begin{equation}\label{firstsumSK2}
-H_{N+n}(\vsi)
=-H_{N+n}(\bolds{\rho}) + \sum_{i\leq n}\eps_i Z_i(\bolds{\rho}),
\end{equation}
where Gaussian processes $Z_i(\vrho)$ have covariance $ \xi'(R(
\bolds{\rho}^{1},\bolds{\rho}^{2}))$.
We can replace the perturbation term $-H_{N+n}^p(\vsi)$ by
%
%e3.5 ###
\begin{equation}\label{HampertSKcav}
-H_{N}^{p}(\vrho) = \sum_{k\leq\pi(c_{N})}\log\ch  G_{\xi',k}(\vrho)
+\sum_{k\leq\pi'(c_{N})} G_{\theta,k}(\vrho)
\end{equation}
without affecting $\mathbb{E}\la h \ra_{N+n}$ asymptotically, since
by Lemma~\ref{LemShiftSK} we can slightly
modify the Poisson number of terms using that $|c_{N+n}-c_N|\to0$
and then replace $G_{\xi',i}(\vsi)$ and $G_{\theta,i}(\vsi)$ by
$G_{\xi',i}(\vrho)$ and $G_{\theta,i}(\vrho)$ by interpolation
using that $c_N=o(N)$.
If $\la\cdot\ra_N'$ denotes the Gibbs average on $\Sigma_N$ with
respect to the Hamiltonian
%
%e3.6 ###
\begin{equation}\label{HNprime}
-H_N'(\bolds{\rho})
=-H_{N+n}(\bolds{\rho})-H_{N}^{p}(\vrho),
\end{equation}
then each factor in
\[
\la h\ra_{N+n}=\prod_{j\leq q} \la h_j \ra_{N+n}
=\prod_{j\leq q} \biggl\la\prod_{i\in C_j}\sigma_i \biggr\ra_{N+n}
=\prod_{j\leq q} \biggl\la\prod_{i\in C_j}\rho_i \biggr\ra_{N+n}
\]
(in the last equality we used that for large $N$ all sets $C_j$ will be
on the first $N$ coordinates)
can be written as
\[
\la h_j \ra_{N+n} = \frac{ \la\prod_{i\in C_j}\rho_i
\Av_\eps\exp\sum_{i\leq n} \eps_i Z_i(\vrho) \ra_N'}{ \la
\Av_\eps\exp\sum_{i\leq n} \eps_i Z_i(\vrho) \ra_N'} = \biggl\la\prod_{i\in
C_j}\rho_i \biggr\ra_N'',
\]
where $\la\cdot\ra_N''$ is the Gibbs average on $\Sigma_N$
corresponding to
the Hamiltonian
\[
-H_N''(\bolds{\rho}) = -H_N'(\bolds{\rho}) + \sum_{i\leq n} \log\ch
Z_i(\vrho).
\]
Thus, $\mathbb{E}\la h\ra_{N+n} = \mathbb{E}\la h\ra_N''$.
Since $Z_i(\vrho)$ are independent copies of $G_{\xi'}(\vrho)$, in
distribution
\[
-H_N''(\vrho) =-H_{N+n}(\vrho)-H_{N}^{p,1}(\vrho),
\]
where
%
%e3.7 ###
\begin{equation}\label{HampertSKcav2}
-H_{N}^{p,1}(\vrho) = \sum_{k\leq\pi(c_{N})+n}\log\ch
G_{\xi',k}(\vrho) +\sum_{k\leq\pi'(c_{N})} G_{\theta,k}(\vrho).\vadjust{\goodbreak}
\end{equation}
Let us now consider $\mathbb{E}\la h\ra_N$. It is easy to check that,
in distribution, the Hamiltonian
$H_N(\vrho)$ can be related to the Hamiltonian $H_{N+n}(\vrho)$ in
(\ref{HamSKcav}) by
%
%e3.8 ###
\begin{equation}\label{secondsumSK}
-H_N(\vrho) = -H_{N+n}(\vrho) + \sum_{i\leq n} Y_i(\vrho),
\end{equation}
where $(Y_i(\vrho))$ are independent Gaussian processes with covariance
\[
\mathbb{E}Y_i(\bolds{\rho}^{1}) Y_i(\bolds{\rho}^{2})
= \theta(R(\bolds{\rho}^{1},\bolds{\rho}^{2})) + o_N(1).
\]
Again, without affecting $\mathbb{E}\la h\ra_N$ asymptotically, one
can assume that the covariance of $Y_i(\vrho)$ is exactly $\theta
(R(\vrho^1,\vrho^2))$ which means that they are independent
copies of $G_{\theta}(\vrho)$. Therefore, we can assume that $\mathbb
{E}\la h\ra_N$ is taken with
respect to the Hamiltonian
\[
-H_N'''(\vrho) =-H_{N+n}(\vrho)-H_{N}^{p,2}(\vrho),
\]
where
%
%e3.9 ###
\begin{equation}\label{HampertSKcav3}
-H_{N}^{p,2}(\vrho) = \sum_{k\leq\pi(c_{N})}\log\ch  G_{\xi',k}(\vrho)
+\sum_{k\leq\pi'(c_{N})+n} G_{\theta,k}(\vrho).
\end{equation}
Lemma~\ref{LemShiftSK} then implies that both perturbation terms (\ref
{HampertSKcav2})
and (\ref{HampertSKcav3}) can be replaced by the original perturbation
term (\ref{HampertSK})
without affecting $\mathbb{E}\la h\ra_N''$ and $\mathbb{E}\la h\ra
_N$ asymptotically and this finishes the proof.
\end{pf}

%s3.2 ###
\subsection{Lower bound} \label{SecLBSK}

\begin{lemma}\label{LemLBSK}
There exists $\mu\in\M$ such that
$
\lim_{N\to\infty} F_N \geq\P(\mu).
$
\end{lemma}
\begin{pf}
We again use (\ref{Zdiff}). Suppose that this lower limit is achieved
over subsequence $(N_k)$
and let $\mu\in\M$ be a limit of $(\mu_N)$ over some subsubsequence
of $(N_k)$.
Let $Z_N'$ and $\la\cdot\ra$ be the partition function and the Gibbs
average on $\Sigma_N$
corresponding to the Hamiltonian $H_{N}'$ defined in (\ref{HNprime}),
and let us compute
the limit of
\[
\mathbb{E}\log\frac{Z_{N+1}}{Z_N'}-\mathbb{E}\log\frac{Z_{N}}{Z_N'}
\]
along the above subsubsequence. Using (\ref{firstsumSK2}) and (\ref
{secondsumSK}) for $n=1$
and the fact that, as in (\ref{HampertSKcav}), the perturbation
Hamiltonian $H_{N+1}^p(\vsi)$ in
$Z_{N+1}$ can be replaced by~$H_N^p(\vrho)$, the above limit is equal
to the limit of
\[
\log2 + \mathbb{E}\log\la\ch  G_{\xi'}(\vrho)\ra- \mathbb
{E}\log\la\exp G_{\theta}(\vrho)\ra.
\]
It remains to show that
%
%e3.10 ###
\begin{equation}\label{limAS1}
\lim_{N\to\infty} \mathbb{E}\log\la\ch  G_{\xi'}(\vrho)\ra=
\mathbb{E}\log\mathbb{E}' \ch  G_{\xi'}(\bar{\sigma}_\mu
(w,u,\cdot))
\end{equation}
and
%
%e3.11 ###
\begin{equation}\label{limAS2}
\lim_{N\to\infty} \mathbb{E}\log\la\exp G_{\theta}(\vrho)\ra=
\mathbb{E}\log\mathbb{E}' \exp G_{\theta}(\bar{\sigma}_\mu
(w,u,\cdot)),\vadjust{\goodbreak}
\end{equation}
where for simplicity of notations we will write limits for $N\to\infty
$ rather than over the above
subsubsequence. The proof of this is identical to Talagrand's proof of the
Baffioni--Rosati theorem in~\cite{SG2}. First of all, if $\mathbb
{E}_g$ denotes the expectation
in the randomness of $G_{\xi'}(\vrho)$ conditionally on the
randomness in $\la\cdot\ra$,
then standard Gaussian concentration implies that (see, e.g., Lemma~3
in~\cite{Posit})
\[
\mathbb{P}_g\bigl(|{\log}\la\ch  G_{\xi'}(\vrho)\ra-
\mathbb{E}_g {\log}\la\ch  G_{\xi'}(\vrho)\ra|
\geq A\bigr)\leq e^{-cA^2}
\]
for some small enough constant $c$, and since
\[
0
\leq
\mathbb{E}_g \log\la\ch  G_{\xi'}(\vrho)\ra
\leq
\log\la\mathbb{E}_g \ch  G_{\xi'}(\vrho)\ra\leq\xi'(1)/2
\]
for large enough $A>0$, we get
%
%e3.12 ###
\begin{equation}\label{concch}
\mathbb{P}\bigl(|{\log}\la\ch  G_{\xi'}(\vrho)\ra | \geq A \bigr)
\leq e^{-cA^2}.
\end{equation}
Therefore, if we denote
$
\log_A x = \max(-A,\min(\log x, A)),
$
then for large\break enough~$A$,
%
%e3.13 ###
\begin{equation}\label{partA}
| \mathbb{E}\log\la\ch  G_{\xi'}(\vrho)\ra - \mathbb{E}\log_A \la\ch
G_{\xi'}(\vrho)\ra | \leq e^{-cA^2}.
\end{equation}
Next, if we define $\ch_A  x = \min(\ch  x,\ch  A)$, then using that
\[
|{\log_A x} - \log_A y| \leq e^A |x-y|
\quad\mbox{and}\quad   | {\ch  x}
- \ch_A  x | \leq\ch  x  I(|x| \geq A)
\]
we can write
\begin{eqnarray*}
| \mathbb{E}\log_A \la\ch  G_{\xi'}(\vrho)\ra - \mathbb{E}\log_A
\la\ch_A  G_{\xi'}(\vrho)\ra | &\leq& e^A \mathbb{E}\la | \ch
G_{\xi'}(\vrho) - \ch_A  G_{\xi'}(\vrho) | \ra
\\
&\leq&
e^A
\mathbb{E}\bigl\la
\ch  G_{\xi'}(\vrho)
I\bigl(|G_{\xi'}(\vrho)| \geq A\bigr)
\bigr\ra.
\end{eqnarray*}
By H\"older's inequality we can bound this by
\[
e^A ( \mathbb{E}\la \mathbb{E}_g \ch^2 G_{\xi'}(\vrho) \ra )^{1/2} \bigl(
\mathbb{E}\bigl\la \mathbb{P}_g \bigl(|G_{\xi'}(\vrho)| \geq A\bigr)
\bigr\ra \bigr)^{1/2} \leq
e^{-cA^2}
\]
for large enough $A$ since $\mathbb{P}_g (|G_{\xi'}(\vrho)| \geq
A)\leq e^{-cA^2}$.
Combining with (\ref{partA}) proves that
%
%e3.14 ###
\begin{equation}\label{partA2}
|\mathbb{E}\log\la\ch  G_{\xi'}(\vrho)\ra -\mathbb{E}\log_A \la\ch_A
G_{\xi'}(\vrho)\ra | \leq e^{-cA^2}.
\end{equation}
Approximating logarithm by polynomials on the interval $[e^{-A},e^A]$
we can approximate
$\mathbb{E}\log_A \la\ch_A  G_{\xi'}(\vrho)\ra$ by some linear
combinations of the moments
\[
\mathbb{E}\la\ch_A  G_{\xi'}(\vrho)\ra^q = \mathbb{E}\biggl\la\prod_{l\leq q}
\ch_A  G_{\xi'}(\vrho^l) \biggr\ra = \mathbb{E}\biggl\la\mathbb{E}_g \prod_{l\leq
q} \ch_A  G_{\xi '}(\vrho^l) \biggr\ra
\]
for $q\geq1$. Since
%
%e3.15 ###
\begin{equation}\label{Fdef}
\mathbb{E}_g \prod_{l\leq q} \ch_A  G_{\xi'}(\vrho^l) =
F((R_{l,l'})_{l,l'\leq q})
\end{equation}
for some continuous bounded function $F$ of the overlaps
$(R_{l,l'})_{l,l'\leq q}$,
Lemma~\ref{LemP1SK} implies that
\[
\lim_{N\to\infty} \mathbb{E}\biggl\la\mathbb{E}_g \prod_{l\leq q} \ch _A
G_{\xi'}(\vrho^l)\biggr\ra = \mathbb{E}F((R_{l,l'}^\infty)_{l,l'\leq q}).\vadjust{\goodbreak}
\]
Let us rewrite the right-hand side in terms of the process $G_{\xi'}$
in (\ref{Gxi}).
Recall the definition of the processes in (\ref{Gxi}) and (\ref{Gtheta}).
If $\mathbb{E}_G$ is the expectation in the Gaussian randomness of
these processes, then
the definition of the function $F$ in (\ref{Fdef}) implies that
\[
\mathbb{E}F((R_{l,l'}^\infty)_{l,l'\leq q}) = \mathbb{E}\mathbb{E}_G
\prod_{l\leq q} \ch_A  G_{\xi'}(\bar {\sigma}_\mu(w,u_l,\cdot)) =
\mathbb{E}(\mathbb{E}' \ch_A  G_{\xi'}(\bar{\sigma}_\mu (w,u,\cdot)))^q
\]
and, therefore,
\[
\lim_{N\to\infty} \mathbb{E}\log_A \la\ch_A  G_{\xi'}(\vrho )\ra =
\mathbb{E}\log_A \mathbb{E}' \ch_A  G_{\xi'}(\bar{\sigma}_\mu
(w,u,\cdot)).
\]
[Notice that this approximation by moments depended on functions of the
overlaps only
which justifies the comment leading to (\ref{Parisimore}).]
One can show similarly to (\ref{partA2}) that
%
%e3.16 ###
\begin{equation} \label{partA3}
| \mathbb{E}\log\mathbb{E}' \ch  G_{\xi'}(\bar{\sigma}_\mu (w,u,\cdot))
- \mathbb{E}\log_A \mathbb{E}' \ch_A  G_{\xi'}(\bar{\sigma}_\mu
(w,u,\cdot)) | \leq e^{-cA^2},\hspace*{-32pt}
\end{equation}
which finishes the proof of (\ref{limAS1}). Equation (\ref{limAS2})
is proved similarly.
\end{pf}

%s3.3 ###
\subsection{Upper bound and free energy}\label{SecUBSK}

Since the perturbation term in (\ref{Hampluspert}) does not affect the limit
of free energy, we will now ignore it and consider free energy $F_N$ defined
for the original unperturbed Hamiltonian (\ref{HamSK}).
\begin{lemma}\label{LemUpperSK}
For any function $\bar{\sigma}\dvtx[0,1]^3\to[-1,+1]$ we have
%
%e3.17 ###
\begin{eqnarray}\label{upperboundSK}
F_N &\leq& \log2 + \frac{1}{N} \mathbb{E}\log \mathbb{E}' \prod_{i\leq N}
\ch G_{\xi',i}(\bar{\sigma}(w,u,\cdot))\nonumber\\[-8pt]\\[-8pt]
&&{} - \frac{1}{N} \mathbb{E}\log
\mathbb{E}' \exp\sum_{i\leq N}
G_{\theta,i}(\bar{\sigma}(w,u,\cdot)).\nonumber
\end{eqnarray}
\end{lemma}
\begin{pf}
This is proved by the Guerra type interpolation as in
\cite{Guerra}.
If, for $t\in[0,1]$, we consider the interpolating Hamiltonian
\begin{eqnarray*}
-H_{N,t}(\vsi) &=& -\sqrt{t} H_{N}(\vsi)
+ \sqrt{1-t} \sum_{i\leq N} \sigma_i G_{\xi',i}(\bar{\sigma
}(w,u,\cdot))\\
&&{}+ \sqrt{t} \sum_{i\leq N} G_{\theta, i}(\bar{\sigma}(w,u,\cdot))
\end{eqnarray*}
and interpolating free energy
\[
\varphi(t) = \frac{1}{N}  \mathbb{E}\log\mathbb{E}' \sum_{\vsi
\in\Sigma_N} \exp(-H_{N,t}(\vsi)) ,
\]
then to prove (\ref{upperboundSK}) it is enough to show that $\varphi
'(t)\leq0$.
This is done by the usual Gaussian integration by parts as in~\cite{Guerra}.
\end{pf}

Before proving invariance properties of Theorem~\ref{ThSCSK} let us
finish the proof
of Theorem~\ref{ThFESK} by showing that if we let
\begin{eqnarray*}
\P_n(\mu) &=& \log2 + \frac{1}{n} \mathbb{E}\log \mathbb{E}' \prod_{i\leq
n} \ch G_{\xi',i}(\bar{\sigma}_\mu(w,u,\cdot))\\
&&{} - \frac{1}{n}
\mathbb{E}\log \mathbb{E}' \exp\sum_{i\leq n}
G_{\theta,i}(\bar{\sigma}_\mu(w,u,\cdot)),
\end{eqnarray*}
then the invariance of Theorem~\ref{ThSCSK} implies the following.
\begin{lemma} \label{LemInvn2SK}
For all $\mu\in\M_{\mathrm{inv}}$, $\P_n(\mu)=\P(\mu)$ for all $n\geq1$.
\end{lemma}
\begin{pf}
The result will follow if we show that for $\bar
{\sigma}=\bar{\sigma}_\mu$
for any $n\geq1$,
%
%e3.18 ###
\begin{equation} \label{ArecSK}
\mathbb{E}\log\frac{\mathbb{E}' \prod_{i\leq n+1} \ch  G_{\xi
',i}(\bar{\sigma}(w,u,\cdot))} {\mathbb{E}' \prod_{i\leq n} \ch
G_{\xi',i}(\bar{\sigma }(w,u,\cdot))} = \mathbb{E}\log\mathbb{E}' \ch
G_{\xi',n+1}(\bar{\sigma }(w,u,\cdot))\hspace*{-32pt}
\end{equation}
and
%
%e3.19 ###
\begin{equation}\label{BrecSK}
\mathbb{E}\log\frac{\mathbb{E}' \exp\sum_{i\leq n+1} G_{\theta
,i}(\bar{\sigma}(w,u,\cdot))} {\mathbb{E}' \exp\sum_{i\leq n}
G_{\theta,i}(\bar{\sigma }(w,u,\cdot))} = \mathbb{E}\log\mathbb{E}'
\exp G_{\theta,n+1}(\bar{\sigma }(w,u,\cdot)).\hspace*{-32pt}
\end{equation}
To prove this we will use invariance properties (\ref{ASCSK}) and
(\ref{preBSCSK}).
Using truncation and Gaussian concentration as in Lemma~\ref{LemLBSK},
to prove
(\ref{ArecSK}) it is enough to show that
\begin{eqnarray*}
&&\mathbb{E}\biggl( \frac{\mathbb{E}' \ch_A   G_{\xi',n+1}(\bar
{\sigma}(w,u,\cdot))\prod_{i\leq n} \ch
G_{\xi',i}(\bar{\sigma}(w,u,\cdot))} {\mathbb{E}' \prod_{i\leq n} \ch
G_{\xi',i}(\bar{\sigma }(w,u,\cdot))} \biggr)^q\\
&&\qquad = \mathbb{E}( \mathbb{E}'
\ch_A   G_{\xi',n+1}(\bar{\sigma }(w,u,\cdot)))^q.
\end{eqnarray*}
Using replicas as in (\ref{Gxil}), the left-hand side can be written as
%
%e3.20 ###
\begin{equation}\label{lhs1}
\mathbb{E} \frac{\mathbb{E}' F \prod_{l\leq q} \prod_{i\leq n} \ch
G_{\xi ',i}(\bar{\sigma}(w,u_l,\cdot))} {(\mathbb{E}' \prod_{i\leq n}
\ch  G_{\xi',i}(\bar{\sigma }(w,u,\cdot)))^q},
\end{equation}
where
\[
F = F((R_{l,l'}^\infty)_{l,l'\leq q}) =
\mathbb{E}_G \prod_{l\leq q} \ch_A   G_{\xi',n+1}(\bar{\sigma
}(w,u_l,\cdot))
\]
is a bounded continuous function of the overlaps defined in (\ref{overlaps}).
Approximating $F$ by polynomials of overlaps and using (\ref{ASCSK})
proves that
(\ref{lhs1}) is equal to
\[
\mathbb{E}F = \mathbb{E}\mathbb{E}_G \prod_{l\leq q} \ch_A
G_{\xi',n+1}(\bar{\sigma}(w,u_l,\cdot))
= \mathbb{E}( \mathbb{E}' \ch_A   G_{\xi',n+1}(\bar{\sigma
}(w,u,\cdot)))^q,
\]
and this finishes the proof of (\ref{ArecSK}). Equation (\ref
{BrecSK}) is proved similarly
using (\ref{preBSCSK}) instead.
\end{pf}

%s3.4 ###
\subsection{Invariance and self-consistency equations}\label{SecSCSK}

\mbox{}

\begin{pf*}{Proofs of Theorems~\ref{ThSCSK} and~\ref{ThSKo}}
Let $h = \prod_{l\leq q} h_l$ where $h_l = \prod_{i\in C_l}\sigma_i^l$.
Consider $\mu\in\M$ which is a limit
of $\mu_N$ over some subsequence $(N_k)$.
By Lem\-ma~\ref{LemShift2SK}, the left-hand side of (\ref{SCSK}) is the
limit of
$\mathbb{E}\la h \ra_{N+n}$ over subsequence $(N_k)$.
The right-hand side of (\ref{SCSK}) will appear as a similar limit
once we rewrite this joint moment
of spins using cavity coordinates. The beginning of the proof will be
identical to the proof of Lemma
\ref{LemShift2SK}, only the spins with coordinates $i\leq n$ will now
play the role of cavity coordinates
instead of spins with coordinates $N+1\leq i\leq N+n$. Let us write
$\vsi=(\bolds{\eps},\vrho)\in\Sigma_{N+n}$
for the first $n$ cavity coordinates $\bolds{\eps}=(\eps_1,\ldots,
\eps_n)$ and the last $N$ coordinates
$\vrho= (\rho_1,\ldots,\rho_N)$. Let us consider sequences of
Gaussian processes $(Z_i(\vrho))$ and $(Y_i(\vrho))$
which are independent copies of $G_{\xi'}(\vrho)$ and $G_{\theta
}(\vrho)$, correspondingly.
First of all, we can replace the perturbation term $-H_{N+n}^p(\vsi)$ with
%
%e3.21 ###
\begin{equation}\label{HampertSKcavI}\quad
-H_{N}^{p}(\vrho) = \sum_{k\leq\pi(c_{N})}\log\ch  G_{\xi',k}(\vrho)
+\sum_{k\leq\pi'(c_{N})} G_{\theta,k}(\vrho) +\sum_{k\leq r}
Y_{k}(\vrho)
\end{equation}
for a fixed $r\geq1$ without affecting $\mathbb{E}\la h
\ra_{N+n}$ asymptotically, since by Lemma~\ref{LemShiftSK} we can slightly
modify the Poisson number of terms, and then we can replace $G_{\xi
',i}(\vsi)$ and $G_{\theta,i}(\vsi)$ with $G_{\xi',i}(\vrho)$ and
$G_{\theta,i}(\vrho)$ by interpolation using that $c_N=o(N)$.
Then, as in (\ref{firstsumSK2}), we can redefine the Hamiltonian
$-H_{N+n}(\vsi)$ by
%
%e3.22 ###
\begin{equation} \label{firstsumSK2I}
-H_{N+n}(\vsi) = -H_{N+n}(\bolds{\rho}) + \sum_{i\leq n}\eps_i
Z_i(\bolds{\rho}),
\end{equation}
where $H_{N+n}(\vrho)$ is defined in (\ref{HamSKcav}).
Let $\la\cdot\ra$ denote the Gibbs average corresponding to the Hamiltonian
%
%e3.23 ###
\begin{equation}\label{here1}\qquad
-H_N'(\vrho) = -H_{N+n}(\bolds{\rho}) + \sum_{k\leq\pi(c_{N})}\log\ch
G_{\xi',k}(\vrho) +\sum_{k\leq\pi'(c_{N})} G_{\theta,k}(\vrho).
\end{equation}
Recalling the relationship (\ref{secondsumSK}) between $H_N(\vrho)$
and $H_{N+n}(\vrho)$,
let us note that Lemma~\ref{LemShiftSK} implies, as in the proof of
Lemma~\ref{LemShift2SK},
that the joint distribution of spins $\mu_N'$ corresponding to the
Hamiltonian (\ref{here1}) converges
to the same limits (over subsequences) as the original sequence\vadjust{\goodbreak}
$\mu_N$.
Let us write the function $h_l(\vsi)$ in terms of $\bolds{\eps}$ and
$\vrho$ as
\[
h_l(\vsi) = \prod_{i\in C_l}\sigma_i
= \prod_{i\in C_l^1}\sigma_i \prod_{i\in C_l^2}\sigma_i
= \prod_{i\in C_l^1}\eps_i \prod_{i\in C_l^2}\rho_i,
\]
where we will abuse the notations and still write $C_l^2$ to denote
the set of coordinates $\rho_i$ corresponding to the original
coordinates $\sigma_{n+i}$.
Then we can write
%
%e3.24 ###
\begin{equation}\label{ratio}
\mathbb{E}\la h \ra_{N+n}
=\mathbb{E} \frac{\prod_{l\leq q} U_{N,l}}{V_N^q},
\end{equation}
where
\[
U_{N,l} = \biggl\la \Av_\eps\prod_{i\in C_l^1} \eps_i \exp\sum_{i\leq
n}\eps_i Z_i(\bolds{\rho}) \prod_{i\in C_l^2} \rho_i \exp\sum_{k\leq r}
Y_{k}(\vrho) \biggr\ra
\]
and
\[
V_N = \biggl\la \Av_\eps\exp\sum_{i\leq n}\eps_i Z_i(\bolds{\rho})
\exp\sum_{k\leq r} Y_{k}(\vrho) \biggr\ra = \la\exp X(\vrho)\ra,
\]
where we introduced
\[
X(\vrho)= \sum_{i\leq n}\log\ch  Z_i(\bolds{\rho})
+ \sum_{k\leq r} Y_{k}(\vrho).
\]
It remains to show that
%
%e3.25 ###
\begin{equation}\label{SCSKagain}
\lim_{N\to\infty} \mathbb{E} \frac{\prod_{l\leq q} U_{N,l}}{V_N^q} =
\mathbb{E}  \frac{\prod_{l\leq q} U_{l}}{V^q},
\end{equation}
where
\[
U_{l} = \mathbb{E}' \Av_\eps\prod_{i\in C_l^1} \eps_i \exp\sum_{i\leq
n}\eps_i G_{\xi',i}(\bar{\sigma}(w,u,\cdot)) \prod_{i\in C_l^2}
\bar{\sigma}_i \exp\sum_{k\leq r} G_{\theta,
k}(\bar{\sigma}(w,u,\cdot))
\]
for $\bar{\sigma}_i = \bar{\sigma}(w,u,v_i)$ and
\[
V = \mathbb{E}' \Av_\eps \exp\sum_{i\leq n}\eps_i
G_{\xi',i}(\bar{\sigma }(w,u,\cdot)) \exp\sum_{k\leq r} G_{\theta,
k}(\bar{\sigma}(w,u,\cdot)),
\]
which is, of course, the same equation as (\ref{SCSK}).
[The proof that (\ref{SS}) implies (\ref{SSSK}) is exactly the same
of the proof of (\ref{SCSKagain}).]
The proof of (\ref{SCSKagain}) is nearly
identical to the proof of (\ref{limAS1}) using truncation and Gaussian
concentration, only
instead of approximating a truncated version of $\log x$ by polynomials
we now need to approximate a
truncated version of $1/x$ by polynomials.
If we denote
\[
Y = \log V_N = \log\la\exp X(\vrho) \ra,
\]
then, as in (\ref{concch}), one can show that for large enough $A>0$
%
%e3.26 ###
\begin{equation}\label{concV}
\mathbb{P}(|Y| \geq A )\leq e^{-cA^2}.\vadjust{\goodbreak}
\end{equation}
For $A>0$ let $(x)_A = \max(-A,\min(x,A))$ so that
\[
|{\exp}(-q x) - \exp(-q(x)_A)| \leq \max(e^{-qA}, \exp(-qx))I(|x|>A).
\]
If we denote $Z = \prod_{l\leq q} U_{N,l}$, then, obviously, $\mathbb
{E}Z^2\leq L$ for some large enough $L>0$
that depends on $q,n,r$ and function $\xi$, and (\ref{concV}) implies that
%
%e3.27 ###
\begin{eqnarray}\label{part1i}
&&|\mathbb{E}Z \exp(-q Y) - \mathbb{E}Z \exp(-q(Y)_A)|
\nonumber\\[-8pt]\\[-8pt]
&&\qquad\leq \mathbb{E}|Z|
\max(e^{-qA}, \exp(-qY))I(|Y|>A)
\leq e^{-c A^2}\nonumber
\end{eqnarray}
for large enough $A$. Next, let $\exp_A x = \max(e^{-A},\min(\exp x,
e^A))$,
and let $Y' = \log\la\exp_A X(\vrho) \ra$. Since for all $x, y\in
\Reals$
\[
|{\exp}(-q(x)_A) - \exp(-q(y)_A)| \leq qe^{(q+1)A}|{\exp x} - \exp y|,
\]
we get
\[
|{\exp}(-q(Y)_A) - \exp(-q(Y')_A)| \leq qe^{(q+1)A} \la|{\exp X}(\vrho) -
\exp_A X(\vrho)| \ra.
\]
Next, since for all $x\in\Reals$
\[
|{\exp x} - \exp_A x| \leq\max(e^{-A}, \exp x) I(|x|\geq A),
\]
we obtain the following bound:
\[
|{\exp}(-q(Y)_A - \exp(-q(Y')_A)| \leq qe^{(q+1)A} \bigl\la\max(e^{-A}, \exp
X(\vrho))I\bigl(|X(\vrho)|\geq A\bigr) \bigr\ra.
\]
It is easy to see that $\mathbb{P}(|X(\vrho)| \geq A )\leq e^{-cA^2}$
for large enough $A$, and using H\"older's inequality,
\begin{eqnarray*}
&&|\mathbb{E}Z \exp(-q(Y)_A) - \mathbb{E}Z \exp(-q(Y')_A)|
\\
&&\qquad\leq
qe^{(q+1)A} (\mathbb{E}Z^2)^{1/2}
(\mathbb{E}\la
\max(e^{-4A}, \exp4X(\vrho))
\ra)^{1/4}
\\
&&\qquad\quad{}
\times
\bigl(\mathbb{E}\bigl\la
I\bigl(|X(\vrho)|\geq A\bigr)
\bigr\ra\bigr)^{1/4}
\\
&&\qquad\leq e^{-cA^2}.
\end{eqnarray*}
Combining this with (\ref{part1i}) we prove that
\[
|\mathbb{E}Z \exp(-q Y) - \mathbb{E}Z \exp(-q(Y')_A)| \leq e^{-c A^2}
\]
for large enough $A$. We can now approximate
$\exp(-q(Y')_A) = \la\exp_A X(\vrho) \ra^{-q} $
uniformly by polynomials of $\la\exp_A X(\vrho) \ra$, and therefore
$\mathbb{E}Z \exp(-q(Y')_A)$ can be approximated by a linear
combination of terms
%
%e3.28 ###
\begin{equation}\label{here2}
\mathbb{E}\prod_{l\leq q} U_{N,l} \la\exp_A X(\vrho)
\ra^s.
\end{equation}
If we write the product of the Gibbs averages using replicas and take
expectation with respect
to the Gaussian processes $(X_i(\vrho))$ and $(Y_i(\vrho))$ inside
the Gibbs average,
we will get a Gibbs average of some bounded continuous function of
finitely many overlaps
in addition to the spin terms $\prod_{i\in C_l^2}\rho_i$ that appear
in the definition of $U_{N,l}$.
Observe that if, from the beginning, we chose $m=n$, then factors
$\prod_{i\in C_l^2}\rho_i$
would not be present, which means that the linear combination of (\ref
{here2}) gives an
approximation of $\mathbb{E}\la h\ra_{N+n}$ (and thus $\mathbb{E}\la
h\ra_N$) by the annealed Gibbs average
of some functions of overlaps only. In particular, this proves Theorem
\ref{ThSKo}.
In the general case, Lemma~\ref{LemP1SK} implies that (\ref{here2})
converges to
$\mathbb{E}\prod_{l\leq q} U_l (V_A)^s$ where
\[
V_A
=
\mathbb{E}'
\exp_A \biggl(\sum_{i\leq n}\log\ch  G_{\xi',i}(\bar{\sigma
}(w,u,\cdot))
+ \sum_{k\leq r} G_{\theta, k}(\bar{\sigma}(w,u,\cdot))\biggr).
\]
Since the same truncation and approximation arguments can be carried
out in parallel for the right-hand side
of (\ref{SCSKagain}), this proves (\ref{SCSKagain}) and finishes the
proof of Theorem~\ref{ThSCSK}.
\end{pf*}
\begin{pf*}{Proof of Theorem~\ref{GGlim}}
Using Gaussian integration by parts and invariance in (\ref{SSSK}),
(\ref{GSS}) can be rewritten as
\[
1 - t^2\bigl(\mathbb{E}(R_{1,2}^\infty)^{2p} - 2\mathbb
{E}(R_{1,2}^\infty)^p (R_{1,3}^\infty)^p
+( \mathbb{E}(R_{1,2}^\infty)^p)^2\bigr),
\]
and the second term disappears whenever the Ghirlanda--Guerra
identities (\ref{GGplim}) hold. On the other hand, if (\ref{GSS}) is
uniformly bounded for all
$t> 0$, then for any bounded continuous function $F$ of multi-overlaps
(\ref{multioverlap}) on $n$ replicas,
%
%e3.29 ###
\begin{eqnarray}\label{GSSF}
&&\mathbb{E}\frac{F G_{p}(\bar{\sigma}(w,u_1,\cdot)) \exp t \sum
_{l\leq n} G_{p}(\bar{\sigma}(w,u_l,\cdot))}
{(\mathbb{E}' \exp t G_{p}(\bar{\sigma}(w,u,\cdot)))^n}
\nonumber\\[-8pt]\\[-8pt]
&&\qquad{}-
\mathbb{E}F
\mathbb{E}\frac{G_{p}(\bar{\sigma}(w,u,\cdot)) \exp t G_{p}(\bar
{\sigma}(w,u,\cdot))}
{\mathbb{E}' \exp t G_{p}(\bar{\sigma}(w,u,\cdot))}\nonumber
\end{eqnarray}
is also uniformly bounded by invariance in (\ref{SSSK}) and H\"older's
inequality.
Using Gaussian integration by parts and invariance in (\ref{SSSK}),
this is equal to
%
%e3.30 ###
\begin{equation}
t\Biggl( \sum_{l=2}^{n} \mathbb{E}F (R_{1,l}^\infty)^p -n
\mathbb{E}F(R_{1,n+1}^\infty)^p + \mathbb{E}F  \mathbb
{E}(R_{1,2}^\infty)^p \Biggr),
\end{equation}
which can be bounded only if (\ref{GGplim}) holds.
\end{pf*}

\section*{Acknowledgments}

The author would like to thank Tim Austin for motivating this work and
Michel Talagrand and anonymous referee for asking good questions and
making many suggestions that helped improve the paper.

%suskaldyti doi

% imsref loaded by lrinkeviciute, 2011-11-09 16:05:12
% imsref loaded by lrinkeviciute, 2011-11-09 16:08:06

\printaddresses

\end{document}